\input ./style/arxiv-general.cfg
\documentclass[aop,MSNbibl,dvips]{arximspdf}
\makeatletter
   \@ifpackageloaded{graphicx}{}{\usepackage{graphicx}}
\makeatother
\doi{10.1214/14-AOP947} %kopijuoti is PTS
\volume{43}
\issue{5}
\pubyear{2015}
\firstpage{2701}
\lastpage{2728}
\docsubty{FLA}

\makeatletter
\def\mid{|}
\newtheorem{theorem}{Theorem}
\newproclaim{definition}{Definition}
\newtheorem{proposition}[theorem]{Proposition}
\newtheorem{lemma}[theorem]{Lemma}

\newproclaim{remark}{Remark}%\theoremstyle{definition}
\newproclaim{remarks}{Remarks}

\def\ve{\varepsilon}
\def\wt{\widetilde}
\def\wh{\widehat}
\def\bR{\mathbf{R}}
\def\la{\longrightarrow}
\def\T{\mathbb{T}}
\def\R{\mathbb{R}}
\def\P{\mathbb{P}}
\def\N{\mathbb{N}}
\def\Z{\mathbb{Z}}
\def\da{\downarrow}
\def\build#1_#2^#3{\mathrel{\mathop{#1}\limits_{#2}^{#3}}}
\newcommand{\eqref}[1]{(\ref{#1})}
\makeatother

\begin{document}
\begin{frontmatter}

\title{The range of tree-indexed random walk in~low~dimensions}
\runtitle{The range of tree-indexed random walk}

\begin{aug}
\author[A]{\fnms{Jean-Fran\c cois}~\snm{Le Gall}\corref{}\ead[label=e1]{jean-francois.legall@math.u-psud.fr}}
\and
\author[A]{\fnms{Shen}~\snm{Lin}\ead[label=e2]{shen.lin.math@gmail.com}}
%
% \thankstext{t2}{Footnote to the first author with the `thankstext'
%command.}
%
\runauthor{J.-F. Le Gall and S. Lin}
\affiliation{Universit\'e Paris-Sud}
\address[A]{D\'epartment de math\'ematiques\\
Universit\'e Paris-Sud\\
91405 Orsay\\
France\\
\printead{e1}\\
\phantom{E-mail:\ }\printead*{e2}}
%\author{\fnms{}~\snm{}\corref{}}
%\and
%\author{\fnms{}~\snm{}}
%\runauthor{}
%\affiliation{}
%\dedicated{}
%\address{} %adresu isvedimo komanda gale!
%\address{}
\end{aug}

% HISTORY:
\received{\smonth{1} \syear{2014}}
%\revised{\smonth{} \syear{}}
%\accepted{\smonth{} \syear{}}

% ABSTRACT
%
\begin{abstract}
We study the range $R_n$ of a random walk on the $d$-dimensional
lattice $\Z^d$ indexed by a random tree with $n$ vertices. Under the
assumption that the random walk is centered and has finite fourth
moments, we prove in dimension $d\leq3$ that $n^{-d/4}R_n$ converges
in distribution to the Lebesgue measure of the support of the
integrated super-Brownian excursion (ISE). An auxiliary result shows
that the suitably rescaled local times of the tree-indexed random walk
converge in distribution to the density process of ISE. We obtain
similar results for the range of critical branching random walk in $\Z
^d$, $d\leq3$. As an intermediate estimate, we get exact asymptotics
for the probability that a critical branching random walk starting with
a single particle at the origin hits a distant point. The results of
the present article complement those derived in higher dimensions in
our earlier work.
\end{abstract}

% KEYWORDS

%
%\begin{keyword}[class=AMS]
%\kwd[Primary ]{}
%\kwd{}
%\kwd[; secondary ]{}
%\end{keyword}
%\begin{keyword}
%\kwd{}
%\end{keyword}
\end{frontmatter}

%s1 #&#
\section{Introduction}\label{sec1}

In the present paper, we continue our study of asymptotics for the
number of distinct sites of the lattice visited by a tree-indexed
random walk.
We consider (discrete) plane trees, which are
rooted ordered trees that can be viewed as describing the genealogy of
a population starting with one ancestor or root, which is denoted by
the symbol $\varnothing$. Given such a tree $\mathcal{T}$ and a
probability measure
$\theta$ on $\Z^d$, we can consider the random
walk with jump distribution $\theta$ indexed by the tree $\mathcal
{T}$. This
means that we assign a
(random) spatial
location $Z_{\mathcal{T}}(u)\in\Z^d$ to every vertex $u$ of
$\mathcal{T}$, in the
following way. First, the
spatial location $Z_{\mathcal{T}}(\varnothing)$ of the root is the
origin of $\Z^d$.
Then we assign
independently to every
edge $e$ of the tree $\mathcal{T}$ a random variable $Y_e$ distributed
according
to $\theta$, and we let
the spatial location $Z_{\mathcal{T}}(u)$ of the vertex $u$ be the sum
of the
quantities $Y_e$ over all edges $e$ belonging to the simple path
from $\varnothing$ to $u$ in the tree. The number of distinct spatial locations
is called the range of the tree-indexed random walk $Z_{\mathcal{T}}$.

In our previous work \cite{LGL}, we stated the following result. Let
$\theta$ be a probability distribution on $\Z^d$, which is symmetric
with finite support and is not supported on a
strict subgroup of $\Z^d$, and for every
integer $n\geq1$, let $\mathcal{T}^\circ_n$ be a random tree
uniformly distributed
over all plane trees with $n$ vertices. Conditionally given
$\mathcal{T}^\circ_n$, let $Z_{\mathcal{T}^\circ_{n}}$ be a random
walk with jump distribution
$\theta$ indexed by $\mathcal{T}^\circ_n$,
and let $R_n$ stand for the range of $Z_{\mathcal{T}^\circ_{n}}$. Then:
\begin{itemize}
\item%[$\bullet$]
if $d\geq5$,
\[
\frac{1}{n} R_n \build{\la}_{n\to\infty}^{\mathrm{(P)}}
c_\theta,
\]
where $c_\theta>0$ is a constant depending on $\theta$, and $\build
{\la
}_{}^{\mathrm{(P)}}$ indicates convergence in probability;
\item%[$\bullet$]
if $d=4$,
\[
\frac{\log n}{n} R_n \build{\la}_{n\to\infty}^{L^2} 8
\pi^2 \sigma^4,
\]
where $\sigma^2=(\operatorname{ det} M_\theta)^{1/4}$, with $M_\theta$ denoting
the covariance matrix of $\theta$;
\item%[$\bullet$]
if $d\leq3$,
%
%
%e1 #&#
\begin{equation}
\label{intro-conv} n^{-d/4} R_n \build{\la}_{n\to\infty}^{\mathrm{(d)}}
c_\theta\lambda_d\bigl(\operatorname{ supp}(\mathcal{I})\bigr),
\end{equation}
where $c_\theta=2^{d/4}(\operatorname{ det} M_{\theta})^{1/2}$ is a constant
depending on $\theta$, and $ \lambda_d(\operatorname{ supp}(\mathcal{I}))$ stands
for the
Lebesgue measure of the support of the random measure on $\R^{d}$ known
as Integrated Super-Brownian Excursion or
ISE (see Section~\ref{sec2.3} below for a definition of ISE in terms
of the
Brownian snake, and note that our
normalization is slightly different from the one in \cite{Al}).
\end{itemize}

Only the cases $d\geq5$ and $d=4$ were proved in \cite{LGL}, in fact
in a greater generality than
stated above, especially when $d\geq5$. In the present work, we
concentrate on the case $d\leq3$ and we
prove a general version of the convergence \eqref{intro-conv}, where
instead of considering a uniformly
distributed plane tree with $n$ vertices we deal with a Galton--Watson
tree with offspring distribution $\mu$
conditioned to have $n$ vertices.

Let us specify the assumptions that will be in force throughout this
article. We always assume that $d\leq3$ and:
\begin{itemize}
\item%[$\bullet$]
$\mu$ is a nondegenerate critical offspring
distribution on $\Z_+$, such that, for some $\lambda>0$,
\[
\sum_{k=0}^\infty e^{\lambda k} \mu(k) <
\infty,
\]
and we set $\rho:=(\operatorname{ var} \mu)^{1/2}>0$;
\item%[$\bullet$]
$\theta$ is a probability measure on $\Z^d$, which is
not supported on a strict subgroup of $\Z^d$; $\theta$ is such that
%
%
%e2 #&#
\begin{equation}
\label{hypoJM} \lim_{r\to+\infty} r^4 \theta\bigl(\bigl\{x
\in\Z^d\dvtx |x|>r\bigr\}\bigr)=0,
\end{equation}
and $\theta$ has zero mean; we set $\sigma:=(\operatorname{ det} M_\theta
)^{1/2d}>0$, where
$M_\theta$ denotes the covariance matrix of $\theta$.
\end{itemize}
Note that \eqref{hypoJM} holds if $\theta$ has finite fourth moments.

For every $n\geq1$ such that this makes sense, let $\mathcal{T}_n$ be a
Galton--Watson tree with offspring
distribution $\mu$ conditioned to have $n$ vertices. Note that the case
when $\mathcal{T}_n$ is uniformly distributed over
plane trees with $n$ vertices is recovered when $\mu$ is the geometric
distribution with
parameter $1/2$ (see, e.g., Section~2.2 in \cite{LGM}). Let
$Z_{\mathcal{T}_n}$
denote the random walk with jump distribution
$\theta$ indexed by $\mathcal{T}_n$, and let $R_n$
be the range of $Z_{\mathcal{T}_n}$. Theorem~\ref{conv-range} below
shows that
the convergence (\ref{intro-conv}) holds, provided that $c_\theta$
is replaced by the constant $2^{d/2}\sigma^d\rho^{-d/2}$.

An interesting auxiliary result is an invariance principle for ``local
times'' of our tree-indexed
random walk. For every $a\in\Z^d$, let
\[
L_n(a)= \sum_{u\in\mathcal{T}_n} \mathbf{1}_{\{Z_{\mathcal
{T}_n}(u)=a\}}
\]
be the number of visits of $a$ by the tree-indexed random walk
$Z_{\mathcal{T}
_n}$. For $x=(x_1,\ldots,x_d)\in\R^d$,
set $\lfloor x\rfloor:=(\lfloor x_1\rfloor,\ldots,\lfloor x_d\rfloor)$.
Then Theorem~\ref{convLT}
shows that the process
\[
\bigl(n^{{d}/{4} -1} L_n\bigl(\bigl\lfloor n^{1/4}x \bigr
\rfloor\bigr) \bigr)_{x\in\R
^d\setminus\{0\}}
\]
converges as $n\to\infty$, in the sense of weak convergence of
finite-dimensional marginals,
to the density process of ISE (up to scaling constants and a linear
transformation of the variable $x$). Notice that the latter density
process exists because
$d\leq3$, by results due to Sugitani \cite{Sug}. In dimension $d=1$,
this invariance principle
has been obtained earlier in a stronger (functional) form by
Bousquet-M\'elou and Janson \cite{BMJ}, Theorem~3.6, in a particular
case, and then
by Devroye and Janson~\cite{DJ}, Theorem~1.1, in a more general setting.
Such a strengthening might also be possible when $d=2$ or $3$, but we
have chosen not
to investigate this question here as it is not relevant to our main
applications. In dimensions $2$ and $3$,
Lalley and Zheng \cite{LZ2}, Theorem~1, also give a closely related
result for local times of critical branching random walk
in the case of a Poisson offspring distribution and for a particular
choice of~$\theta$.

Our tree-indexed random walk can be viewed as a branching random walk
starting with
a single initial particle and conditioned to have a fixed total
progeny. Therefore, it is not
surprising that our main results have analogs for branching random
walks, as
it was already the case in dimension $d\geq4$ (see Propositions~20 and~21
in \cite{LGL}). For every integer
$p\geq1$, consider a (discrete time) branching random walk starting
initially with
$p$ particles located at the origin of $\Z^d$, such that the offspring number
of each particle is distributed according to $\mu$, and each newly born
particle jumps
from the location of its parent according to the jump
distribution $\theta$. Let $\mathcal{V}^{[p]}$
stand for the set of all sites of $\Z^d$ visited by this branching
random walk. Then Theorem~\ref{rangeBRW}
shows that, similarly as in \eqref{intro-conv}, the asymptotic
distribution of $p^{-d/2}\#\mathcal{V}^{[p]}$ is the Lebesgue
measure of the range of a super-Brownian motion starting from $\delta
_0$ (note again
that this Lebesgue measure is positive because $d\leq3$, see \cite
{DIP} or \cite{Sug}). In a related direction, we mention the article of
Lalley and Zheng~\cite{LZ},
which gives estimates for the number of occupied sites \textit{at a given
time} by a critical nearest neighbor
branching random walk in $\Z^d$.

Our proof of Theorem~\ref{rangeBRW} depends on an asymptotic estimate
for the hitting probability
of a distant point by branching random walk, which seems to be new and
of independent interest. To be specific, consider
the set $\mathcal{V}^{[1]}$ of all sites visited by the branching
random walk starting with a single
particle at the origin. Consider for simplicity the isotropic case
where $M_\theta=\sigma^2 \operatorname{ Id}$, where $\operatorname{ Id}$ is the identity
matrix. Then Theorem~\ref{estim-visit}
shows that
\[
\lim_{|a|\to\infty} |a|^2 P \bigl(a\in\mathcal{V}^{[1]}
\bigr) = \frac
{2(4-d)\sigma^2}{\rho^2}.
\]
See Section~\ref{sec5.1} for a discussion of similar estimates in
higher dimensions.

Not surprisingly, our proofs depend on the known relations between
tree-indexed random walk
(or branching random walk) and the Brownian snake (or super-Brownian
motion). In particular,
we make extensive use of a result of Janson and Marckert \cite{JM}
showing that the ``discrete
snake'' coding our tree-indexed random walk $Z_{\mathcal{T}_n}$
converges in
distribution in a strong (functional) sense
to the Brownian snake driven by a normalized Brownian excursion. It
follows from this convergence
that the set of all sites visited by the tree-indexed random walk
converges in distribution
(modulo a suitable rescaling) to the support of ISE, in the sense of
the Hausdorff distance between
compact sets. But, of course, this is not sufficient to derive
asymptotics for the \textit{number}
of visited sites.

Our assumptions on $\mu$ and $\theta$ are similar to those in \cite
{JM}. We have not striven for
the greatest generality, and it is plausible that these assumptions can
be relaxed.
See, in particular, \cite{JM} for a discussion of the necessity of the
existence of exponential
moments for the offspring distribution $\mu$. It might also be possible
to replace our condition \eqref{hypoJM} on $\theta$
by a second moment assumption, but this would require different methods
as the results of \cite{JM} show that the strong convergence of
discrete snakes to the Brownian snake
no longer holds without \eqref{hypoJM}.

The paper is organized as follows. Section~\ref{sec2} presents our
main notation
and gives some preliminary
results about the Brownian snake. Section~\ref{sec3} is devoted to our main
result about the
range of tree-indexed random walk in dimension \mbox{$d\leq3$}. Section~\ref{sec4}
discusses similar
results for branching random walk, and Section~\ref{sec5} presents a few
complements and open questions.

%s2 #&#
\section{Preliminaries on trees and the Brownian snake}\label{sec2}
\label{preli}

%s2.1 #&#
\subsection{Finite trees}\label{sec2.1}
\label{fitree}

We use the standard formalism for plane trees. We set
\[
\mathcal{U}:=\bigcup_{n=0}^\infty
\N^n,
\]
where $\N=\{1,2,\ldots\}$ and $\N^0=\{\varnothing\}$. If
$u=(u_1,\ldots,u_n)\in\mathcal{U}$,
we set $|u|=n$ [in particular $|\varnothing|=0$].
We write $\prec$ for the lexicographical order on $\mathcal{U}$, so
that $\varnothing\prec1
\prec(1,1) \prec2$, for instance.

%If $u,v\in\mathcal{U}$, $uv$ stands for the
%concatenation of $u$ and $v$. In particular $\varnothing u=u
%\varnothing=u$.
%The genealogical (partial) order $\ll$ is then defined by
%saying that $u \ll v$ if and only if
%$v=uw$ for some $w\in\mathcal{U}$.

A plane tree (or rooted ordered tree) $ \mathcal{T}$ is a finite
subset of
$\mathcal{U}$
such that:
\begin{longlist}[(iii)]
\item[(i)] $\varnothing\in\mathcal{T}$;

\item[(ii)] If $u=(u_1,\ldots,u_n)\in\mathcal{T}\setminus\{
\varnothing\}$ then
$\check u:= (u_1,\ldots,u_{n-1})\in\mathcal{T}$;

\item[(iii)] For every $u=(u_1,\ldots,u_n)\in\mathcal{T}$, there
exists an
integer $k_u(\mathcal{T})\geq0$
such that, for every $j\in\N$, $(u_1,\ldots,u_n,j)\in\mathcal{T}$
if and only if
$1\leq j\leq
k_u(\mathcal{T})$.
\end{longlist}
The notions of a descendant or of an ancestor of a vertex of $\mathcal
{T}$ are
defined in an obvious way. If $u,v\in\mathcal{T}$, we
will write $u\wedge v\in\mathcal{T}$ for the most recent common
ancestor of $u$
and $v$. We denote the set of all planes trees by $\T_f$.

%The notions of a child and a parent of a vertex of $\t$ are defined in
%an obvious way.
%The quantity $k_u(\t)$ in (iii) is the number of children of $u$ in $
%\t$.
%If $u\in\t$, we write $[\t]_u=\{v\in\mathcal{U}:uv\in\t\}$, which
%corresponds to the
%subtree of descendants of $u$ in $\t$. We denote the set of all plane
%trees by $\T_f$.

Let $\mathcal{T}$ be a tree with $p=\#\mathcal{T}$ vertices and let
$\varnothing
=v_0\prec v_1\prec\cdots\prec v_{p-1}$ be
the vertices of $\mathcal{T}$ listed in lexicographical order. We
define the height function $(H_i)_{0\leq i\leq p}$
of $\mathcal{T}$
by setting $H_i=|v_i|$ for every $0\leq i\leq p-1$, and $H_p=0$ by convention.

Recall that we have fixed a probability measure $\mu$ on $\Z_+$
satisfying the assumptions given in Section~\ref{sec1}, and that $\rho
^2=\operatorname{
var} \mu$. The law of the Galton--Watson tree with offspring
distribution $\mu$ is a probability measure on the space $\T_f$, which
is denoted by $\Pi_\mu$ (see, e.g., \cite{LG1}, Section~1).

We will need some information about the law of the total progeny
$\#\mathcal{T}$ under~$\Pi_\mu$.
It is well known (see, e.g., \cite{LG1}, Corollary~1.6) that this law
is the same as the law of the first hitting time of $-1$ by
a random walk on $\Z$ with jump distribution $\nu(k)=\mu(k+1),
k=-1,0,1,\ldots$
started from $0$. Combining this with Kemperman's formula (see,
e.g., \cite{pitman}, page 122)
and using a standard local limit theorem, one gets
%
%
%e3 #&#
\begin{equation}
\label{Kemp1} \lim_{k\to\infty} k^{1/2}
\Pi_\mu(\#\mathcal{T}\geq k)= \frac
{2}{\rho\sqrt{2\pi}}.
\end{equation}
Suppose that $\mu$ is not supported on a strict subgroup of $\Z$, so
that the
random walk with jump distribution $\nu$ is aperiodic. The preceding
asymptotics can then
be strengthened in the form
%
%
%e4 #&#
\begin{equation}
\label{Kemp2} \lim_{k\to\infty} k^{3/2}
\Pi_\mu(\#\mathcal{T}= k)= \frac
{1}{\rho\sqrt{2\pi}}.
\end{equation}

%s2.2 #&#
\subsection{Tree-indexed random walk}\label{sec2.2}
\label{TRW}

A ($d$-dimensional) spatial tree is a pair $(\mathcal{T},(z_u)_{u\in
\mathcal{T}})$
where $\mathcal{T}\in\T_f$ and $z_u\in\Z^d$ for every $u\in
\mathcal{T}$. We let $\T_f^*$
be the set of all spatial trees.

Recall that $\theta$ is a probability measure on $\Z^d$ satisfying the
assumptions listed in the \hyperref[sec1]{Introduction}.
We write
$\Pi^*_{\mu,\theta}$ for the probability distribution on $\T_f^*$ under
which $\mathcal{T}$ is distributed according to $\Pi_\mu$ and,
conditionally on $\mathcal{T}$, the ``spatial locations'' $(z_u)_{u\in
\mathcal{T}}$ are
distributed as random walk indexed by $\mathcal{T}$,
with jump distribution $\theta$, and started from $0$ at the root
$\varnothing$: This means that,
under the probability measure $\Pi^*_{\mu,\theta}$, we have
$z_\varnothing=0$ a.s. and,
conditionally on $\mathcal{T}$, the quantities $(z_u-z_{\check u},
u\in\mathcal{T}
\setminus\{\varnothing\})$ are independent
and distributed according to $\theta$.

%s2.3 #&#
\subsection{The Brownian snake}\label{sec2.3}

We refer to \cite{Zurich} for the basic facts about the Brownian snake
that we will use. The Brownian snake
$(W_s)_{s\geq0}$ is
a Markov process taking values in the space $\mathcal{W}$ of all
($d$-dimensional) stopped paths: Here,
a~stopped path $w$ is just a continuous mapping $w\dvtx [0,\zeta_{(w)}] \la
\R^d$, where the
number $\zeta_{(w)}\geq0$, which depends on $w$, is called the
lifetime of $w$. A stopped path
$w$ with zero lifetime will be identified with its starting point
$w(0)\in\R^d$. The endpoint $w(\zeta_{(w)})$
of a stopped path $w$ is denoted by $\wh w$.

It will be convenient to argue on the canonical
space $C(\R_+, \mathcal{W})$ of all continuous mappings from $\R_+$
into $\mathcal{W}$, and
to let $(W_s)_{s\geq0}$ be the canonical process on this space. We
write $\zeta_s:= \zeta_{(W_s)}$ for the
lifetime of $W_s$. If $x\in\R^d$, the law of the Brownian snake
starting from $x$ is the probability measure $\P_x$
on $C(\R_+, \mathcal{W})$ that is characterized as follows:
\begin{enumerate}[(ii)]
\item[(i)] The distribution of $(\zeta_s)_{s\geq0}$ under $\P_x$ is
the law of a reflected linear Brownian
motion on $\R_+$ started from $0$.
\item[(ii)] We have $W_0=x$, $\P_x$ a.s. Furthermore, under $\P_x$ and
conditionally on $(\zeta_s)_{s\geq0}$,
the process $(W_s)_{s\geq0}$ is (time-inhomogeneous) Markov with
transition kernels specified as follows. If
$0\leq s<s'$,
\begin{itemize}
\item%[$\bullet$]
$W_{s'}(t)=W_s(t)$ for every $0\leq t\leq m_\zeta
(s,s'):= \min\{\zeta_r\dvtx  s\leq r\leq s'\}$;
\item%[$\bullet$]
$(W_{s'}(m_{\zeta}(s,s')+t)-W_{s'}(m_{\zeta
}(s,s')))_{0\leq t\leq\zeta_{s'}-m_\zeta(s,s')}$ is
a standard Brownian motion in $\R^d$ independent of $W_s$.
\end{itemize}
\end{enumerate}
We will refer to the process $(W_s)_{s\geq0}$ under $\P_0$ as the
standard Brownian snake.

We will also be interested in (infinite) excursion
measures of the Brownian snake, which we denote by $\N_x$, $x\in\R^d$.
For every
$x\in\R^d$, the distribution of the process $(W_s)_{s\geq0}$
under $\N_x$ is characterized by properties analogous to (i) and~(ii)
above, with the
only difference that in (i) the law of reflected linear Brownian
motion is replaced by the It\^o measure of positive excursions of linear
Brownian motion, normalized in such a way that $\N_x(\sup\{\zeta
_s\dvtx s\geq0\}>\ve)=(2\ve)^{-1}$, for every $\ve>0$.

We write $\gamma:= \sup\{s\geq0 \dvtx \zeta_s>0\}$, which
corresponds to the duration of the excursion
under $\N_x$. A special role will be played by the probability measures
$\N^{(r)}_x:= \N_x(\cdot\mid\gamma=r)$, which are defined
for every $x\in\R^d$ and every $r>0$.
Under $\N^{(r)}_x$, the ``lifetime process'' $(\zeta_s)_{0\leq s\leq
r}$ is a Brownian excursion with duration $r$. From the analogous
decomposition for the It\^o measure of Brownian excursions, we have
%
%
%e5 #&#
\begin{equation}
\label{decoIto} \N_0 = \int_0^\infty
\frac{\mathrm{d}r}{2\sqrt{2\pi r^3}} \N^{(r)}_0.
\end{equation}

The total occupation measure of the Brownian snake is the finite
measure $\mathcal{Z}$ on $\R^d$ defined under $\N_x$,
or under $\N^{(r)}_x$, by the formula
\[
\langle\mathcal{Z}, \varphi\rangle= \int_0^\gamma
\,\mathrm{d}s\, \varphi(\wh W_s),
\]
for any nonnegative measurable function $\varphi$ on $\R^d$.

Under
$\N^{(1)}_x$, $\mathcal{Z}$ is a random probability measure, which in
the case
$x=0$ is
called ISE for integrated super-Brownian excursion [the measure
$\mathcal{I}$ in \eqref{intro-conv} is
thus distributed as $\mathcal{Z}$ under $\N^{(1)}_0$]. Note that our
normalization of
ISE is slightly different from the one originally proposed by Aldous
\cite{Al}.

The following result will be derived from known properties of
super-Brownian motion via
the connection between the Brownian snake and superprocesses.

%
%pr1 #&#
\begin{proposition}
\label{densityZ}
Both $\N_x$ a.e. and $\N^{(1)}_x$ a.s., the random measure
$\mathcal{Z}$ has a continuous density on $\R^d$, which will
be denoted by $(\ell^y,y\in\R^d)$.
\end{proposition}

\begin{remark*}
When $d=1$, this result, under the measure $\N^{(1)}_0$,
can be found in~\cite{BMJ}, Theorem~2.1.
\end{remark*}

\begin{pf*}{Proof of Proposition \ref{densityZ}} By translation invariance, it is enough to consider the case
$x=0$. We rely on the Brownian snake construction of
super-Brownian motion to deduce the statement of the proposition from Sugitani's
results \cite{Sug}. Let $(W^i)_{i\in I}$ be a Poisson point measure
on $C(\R_+,\mathcal{W})$ with intensity $\N_0$. With every $i\in I$,
we associate the occupation measure $\mathcal{Z}^i$ of $W^i$.
Then Theorem IV.4 in \cite{Zurich} shows that there exists
a super-Brownian motion $(X_t)_{t\geq0}$ with branching mechanism
$\psi(u)=2u^2$ and initial value $X_0=\delta_0$, such that
\[
\int_0^\infty\,\mathrm{d}t\, X_t = \sum
_{i\in I} \mathcal{Z}^i.
\]
As a consequence of \cite{Sug}, Theorems 2 and 3, the random measure
$\int_0^\infty\,\mathrm{d}t\, X_t $ has a.s. a continuous density on $\R
^d\setminus\{0\}$. On the other hand, let
$B(0,\ve)$ denote the closed ball of radius $\ve$ centered at $0$ in
$\R
^d$. Then,
for every $\ve>0$, the event
\[
\mathcal{A}_\ve:= \bigl\{\#\bigl\{i\in I\dvtx \mathcal{Z}^i
\bigl(B(0,\ve)^c\bigr) >0\bigr\} =1\bigr\}
\]
has positive probability (see, e.g., \cite{Zurich}, Proposition V.9).
On the event $\mathcal{A}_\ve$,
write $i_0$ for the unique index in $I$
such that $\mathcal{Z}^{i_0}(B(0,\ve)^c) >0$. Then, still on the
event $\mathcal
{A}_\ve$, the measures $\int_0^\infty\,\mathrm{d}t\, X_t $
and $\mathcal{Z}^{i_0}$ coincide on $B(0,\ve)^c$. The conditional distribution
of $W^{i_0}$ knowing $\mathcal{A}_\ve$ is $\N_0(\cdot\mid\mathcal
{Z}(B(0,\ve
)^c)>0)$, and we
conclude that $\mathcal{Z}$ has a continuous density on $B(0,\ve)^c$,
$\N_0(\cdot\mid \mathcal{Z}(B(0,\ve)^c)>0)$ a.s. As this holds for
any $\ve
>0$, we obtain
that, $\N_0$ a.e., the random measure
$\mathcal{Z}$ has a continuous density on $\R^d\setminus\{0\}$. Via
a scaling argument,
the same property holds $\N^{(1)}_0$ a.s. This argument does not
exclude the possibility that
$\mathcal{Z}$ might have a singularity at $0$, but we can use the rerooting
invariance property (see \cite{Al}, Section~3.2 or \cite{LGW},
Section~2.3) to complete the proof.
According to this property, if under the measure $\N^{(1)}_0$
we pick a random point distributed according to $\mathcal{Z}$ and then
shift $\mathcal{Z}
$ so that this random point becomes the
origin of $\R^d$, the resulting random measure has the same
distribution as $\mathcal{Z}$. Consequently, we obtain that $\N
^{(1)}_0$ a.s.,
$\mathcal{Z}(\mathrm{d}x)$ a.e., the measure $\mathcal{Z}$ has a
continuous density on $\R
^d\setminus\{x\}$. It easily follows that
$\mathcal{Z}$ has a continuous density on $\R^d$, $\N^{(1)}_0$ a.s.,
and by
scaling again the same property
holds under $\N_0$.
\end{pf*}

Let us introduce the random closed set
\[
\mathcal{R}:= \{\wh W_s \dvtx 0\leq s \leq\gamma \}.
\]
Note that, by construction, $\mathcal{Z}$ is supported on $\mathcal
{R}$, and it
follows that,
for every $y\in\R^d\setminus\{x\}$,
%
%
%e6 #&#
\begin{equation}
\label{inclusion-hitting} \bigl\{\ell^y >0\bigr\} \subset\{y\in\mathcal{R}\},\qquad
\N_x\mbox{ a.e. or }\N ^{(1)}_x\mbox{ a.s.}
\end{equation}

%
%pr2 #&#
\begin{proposition}
\label{hitting-point}
For every $y\in\R^d\setminus\{x\}$,
\[
\bigl\{\ell^y >0\bigr\} = \{y\in\mathcal{R}\},\qquad \N_x
\mbox{ a.e. and }\N ^{(1)}_x\mbox{ a.s.}
\]
\end{proposition}

\begin{pf}
Fix $y\in\R^d$, and consider the function
$u(x):=\N_x(\ell^y>0)$, for every $x\in\R^d\setminus\{y\}$. By
simple scaling
and rotational invariance
arguments (see the proof of Proposition V.9(i) in \cite{Zurich} for a
similar argument), we have
\[
u(x)=C_d|x-y|^{-2}
\]
with a certain constant $C_d>0$ depending only on $d$.
On the other hand, an easy application of
the special Markov property \cite{LG0} shows that, for
every $r>0$, and every $x\in B(y,r)^c$, we have
\[
u(x)=\N_x \biggl[1-\exp \biggl(-\int X^{B(y,r)^c}(\mathrm{d}z)
u(z) \biggr) \biggr],
\]
where $X^{B(y,r)^c}$ stands for the exit measure of the Brownian snake
from the
open set $B(y,r)^c$. Theorem V.4 in \cite{Zurich} now shows that
the function $u$ must solve the partial differential
equation $\Delta u= 4 u^2$ in $\R^d\setminus\{y\}$. It
easily follows that $C_d=2 - d/2$.

The preceding line of reasoning
also applies to the function $v(x):= \N_x(y\in\mathcal{R})$
(see \cite{Zurich}, page 91), and shows
that we have $v(x)=(2-d/2)|x-y|^{-2}= u(x)$
for every $x\in\R^d\setminus\{y\}$---note that this formula for $v$
can also be
derived from \cite{DIP}, Theorem~1.3 and the connection between the
Brownian snake and
super-Brownian motion. Recalling \eqref{inclusion-hitting}, this is
enough to conclude
that
%
%
%e7 #&#
\begin{equation}
\label{hitting-tech0} \bigl\{\ell^y >0\bigr\} = \{y\in\mathcal{R}\},\qquad
\N_x\mbox{ a.e.}
\end{equation}
for every $x\in\R^d\setminus\{y\}$.

We now want to obtain that the equality in \eqref{hitting-tech0} also
holds $\N^{(1)}_x$ a.s.
Note that, for every fixed $x$, we could use a scaling argument to get that
$\{\ell^y >0\} = \{y\in\mathcal{R}\}$, $\N^{(1)}_x$ a.s., for
$\lambda
_d$ a.e. $y\in\R^d$, where we recall that
$\lambda_d$ stands for Lebesgue measure on~$\R^d$. In order
to get the more precise assertion of the proposition, we use a
different method.

By translation invariance, we may assume that $x=0$ and we fix $y\in\R
^d\setminus\{0\}$.
We set $T_y:=\inf\{s\geq0 \dvtx \wh W_s=y\}$. Also, for every $s>0$,
we set
\[
\wt\ell^y_s:= \liminf_{\ve\to0} \bigl(
\lambda_d\bigl(B(y,\ve)\bigr) \bigr)^{-1}\int
_0^s \,\mathrm{d}r\, \mathbf{1}_{\{|\wh W_r-y|\leq\ve\}}.
\]
Note that $\wt\ell^y_\gamma= \ell^y$, $\N_0$ a.e. and $\N
^{(1)}_0$ a.e.
We then claim that, for every $s>0$,
%
%
%e8 #&#
\begin{equation}
\label{hitting-tech1} \{T_y\leq s \} = \bigl\{\wt\ell^y_s
>0 \bigr\},\qquad \N_0\mbox{ a.e.}
\end{equation}
The inclusion $\{\wt\ell^y_s >0\}\subset\{T_y\leq s\} $ is obvious.
In order to prove the
reverse inclusion, we argue by contradiction and assume that
\[
\N_0\bigl(T_y\leq s, \wt\ell^y_s
=0\bigr) >0.
\]
Note that $\N_0(T_y=s)=0$ [because $\N_0(\wh W_s=y)=0$], and so we
have also
$\N_0(T_y< s,\wt\ell^y_s =0) >0$. For every $\eta>0$, let
\[
T^{(\eta)}_y:= \inf \bigl\{r\geq T_y \dvtx
\zeta_r\leq(\zeta _{T_y} - \eta)^+ \bigr\}.
\]
Notice that, by the properties of the Brownian snake, the path
$W_{T^{(\eta)}_y}$ is just
$W_{T_y}$ stopped at time $(\zeta_{T_y}-\eta)^+$.

From the strong Markov property at time $T_y$, we easily get that
$T^{(\eta)}_y\da T_y$
as $\eta\da0$, $\N_0$ a.e. on $\{T_y<\infty\}$. Hence, on the event
$\{
T_y<s\}$, we have also $T^{(\eta)}_y<s$
for $\eta$ small enough, $\N_0$ a.e. Therefore, we can find $\eta>0$
such that
\[
\N_0 \bigl(T_y< s,\wt\ell^y_{T^{(\eta)}_y}
=0 \bigr) >0.
\]

However, using the strong
Markov property at time $T^{(\eta)}_y$, and Lemma V.5 and Proposition
V.9(i) in \cite{Zurich}, we immediately
see that, conditionally on the past up to time $T^{(\eta)}_y$, the
event $\{\wh W_r\neq y, \forall r\geq T^{(\eta)}_y\}$
occurs with positive probability. Hence, we get
\[
\N_0\bigl(T_y<s, \wt\ell^y_\gamma=0
\bigr) >0.
\]
Since $\wt\ell^y_\gamma= \ell^y$, this contradicts \eqref
{hitting-tech0}, and this contradiction completes
the proof of our claim \eqref{hitting-tech1}.

Finally, we observe that, for every $s\in(0,1)$, the law of
$(W_r)_{0\leq r\leq s}$ under $\N^{(1)}_0$
is absolutely continuous with respect to the law of the same process
under $\N_0$ (this is
a straightforward consequence of the similar property for the It\^o
excursion measure and the
law of the normalized Brownian excursion; see, e.g., \cite{RY}, Chapter~XII). Hence, \eqref{hitting-tech1} also gives, for every
$s\in(0,1)$,
\[
\{T_y\leq s \} = \bigl\{\wt\ell^y_s >0
\bigr\},\qquad \N ^{(1)}_0\mbox{ a.s.},
\]
and the fact that the equality in \eqref{hitting-tech0} also holds $\N
^{(1)}_0$ a.s. readily follows.
\end{pf}

%s3 #&#
\section{Asymptotics for the range of tree-indexed random walk}\label{sec3}

Throughout this section, we consider only integers $n\geq1$ such that
$\Pi_\mu( \#\mathcal{T}=n)>0$
(and when we let $n\to\infty$, we mean along such values).
For every such integer $n$,
let $(\mathcal{T}_n,(Z^n(u))_{u\in\mathcal{T}_n})$ be distributed
according to
$\Pi^*_{\mu,\theta}(\cdot\mid \#\mathcal{T}=n)$. Then $\mathcal
{T}_n$ is a
Galton--Watson tree with offspring
distribution $\mu$ conditioned to have $n$ vertices, and conditionally
on $\mathcal{T}_n$,
$(Z^n(u))_{u\in\mathcal{T}_n}$ is a random walk with jump distribution
$\theta$
indexed by $\mathcal{T}_n$.

We set, for every $t> 0$ and $x\in\mathbb{R}^d$,
\[
p_t(x):= \frac{1}{(2\pi t)^{d/2} \sqrt{\operatorname{ det} M_\theta}} \exp \biggl( -\frac{x\cdot M_\theta^{-1}x}{2 t} \biggr),
\]
where $x\cdot y$ stands for the usual scalar product in $\R^d$.

For every $a\in\mathbb{Z}^d$, we also set
\[
L_n(a):= \sum_{u\in\mathcal{T}_n}
\mathbf{1}_{\{Z^n(u)=a\}}.
\]

%
%le3 #&#
\begin{lemma}
\label{contiTlocal}
For every $\ve>0$, there exists a constant $C_\ve$ such that, for
every $n$
and every $b\in\mathbb{Z}^d$ with $|b|\geq\ve n^{1/4}$,
\[
E \bigl[\bigl(L_n(b)\bigr)^2 \bigr] \leq
C_\ve n^{2-{d}/{2}}.
\]
Furthermore, for every $x,y\in\mathbb{R}^d\setminus\{0\}$, and for
every choice of the sequences $(x_n)$ and $(y_n)$ in $\mathbb{Z}^d$
such that $n^{-1/4}x_n\la x$ and $n^{-1/4}y_n\la y$ as $n\to\infty$,
we have
\[
\lim_{n\to\infty} n^{{d}/{2}-2} E \bigl[ L_n(x_n)
L_n(y_n) \bigr]=\varphi(x,y),
\]
where
\begin{eqnarray*}
\varphi(x,y)&:=& {\rho^4} \int_{(\mathbb{R}_+)^3}
\,\mathrm{d}r_1 \,\mathrm {d}r_2 \,\mathrm{d}r_3
(r_1+r_2+r_3)e^{-\rho^2(r_1+r_2+r_3)^2/2}\\
&&{}\times \int
_{\mathbb{R}^d} \,\mathrm{d}z\, p_{r_1}(z)p_{r_2}(x-z)p_{r_3}(y-z).
\end{eqnarray*}
The function $\varphi$ is continuous on $(\mathbb{R}^d\setminus\{0\})^2$.
\end{lemma}

\begin{remark*} The function $\varphi$ is in fact continuous on $(\R
^d)^2$. Since we will not need this result, we leave the
proof to the reader.
\end{remark*}

\begin{pf*}{Proof of Lemma \ref{contiTlocal}} We first establish the second assertion of the lemma. We let
$u^n_0,u^n_1,\ldots,u^n_{n-1}$ be the vertices of $\mathcal{T}_n$
listed in lexicographical order. By definition,
\[
L_n(x_n)=\sum_{i=0}^{n-1}
\mathbf{1}_{\{Z^n(u^n_i)=x_n\}},
\]
so that
\[
E \bigl[L_n(x_n)L_n(y_n) \bigr]
=E \Biggl[ \sum_{i=0}^{n-1}\sum
_{j=0}^{n-1} \mathbf{1}_{\{Z^n(u^n_i)=x_n, Z^n(u^n_j)=y_n\}} \Biggr].
\]
Let $H^n$ be the height function of the tree $\mathcal{T}_n$, so that
$H^n_i=|u^n_i|$ for every
$i\in\{0,1,\ldots,n-1\}$. If $i,j\in\{0,1,\ldots,n-1\}$, we also use
the notation $\check H^n_{i,j}=|u^n_i\wedge u^n_j|$ for the generation
of the most
recent common ancestor to $u^n_i$ and $u^n_j$, and note that
%
%
%e9 #&#
\begin{equation}
\label{MRCA} \Bigl|\check H^n_{i,j} - \min
_{i\wedge j\leq k\leq i\vee j} H^n_k \Bigr| \leq1.
\end{equation}

Write $\pi_k=\theta^{*k}$ for the transition kernels of
the random walk with jump distribution $\theta$.
By conditioning with respect to the tree $\mathcal{T}_n$, we get
%
%
%e10 #&#
\begin{eqnarray}
\label{Tlocal1}&& E \bigl[L_n(x_n)L_n(y_n)
\bigr]\nonumber\\
&&\qquad=E \Biggl[ \sum_{i=0}^{n-1}\sum
_{j=0}^{n-1} \sum
_{a\in\mathbb{Z}^d} \pi_{\check H^n_{i,j}}(a) \pi_{H^n_i-\check H^n_{i,j}}(x_n-a)
\pi _{H^n_j-\check H^n_{i,j}}(y_n-a) \Biggr]
\\
&&\qquad= n^2 E \biggl[\int_0^1\int
_0^1 \,\mathrm{d}s \,\mathrm{d}t\, \Phi
^n_{x_n,y_n} \bigl(H^n_{\lfloor ns\rfloor},H^n_{\lfloor nt\rfloor},
\check H^n_{\lfloor ns\rfloor,\lfloor nt\rfloor} \bigr) \biggr],\nonumber
\end{eqnarray}
where we have set, for every integers $k,\ell,m\geq0$ such that
$k\wedge\ell\geq m$,
\[
\Phi^n_{x_n,y_n}(k,\ell,m):= \sum
_{a\in\mathbb{Z}^d}\pi_m(a) \pi _{k-m}(x_n-a)
\pi_{\ell-m}(y_n-a).
\]

In the remaining part of the proof, we assume that $\theta$ is
aperiodic [meaning that the
subgroup generated by $\{k\geq0\dvtx \pi_k(0)>0\}$ is $\Z$]. Only
minor modifications are
needed to treat the general case. We can then use
the local limit theorem, in a form that can be obtained by combining
Theorems 2.3.9 and 2.3.10
in \cite{LL}. There exists a sequence $\delta_n$ converging to $0$ such
that, for every $n\geq1$,
%
%
%e11 #&#
\begin{equation}
\label{LLT} \sup_{a\in\mathbb{Z}^d} \biggl( \biggl(1 +
\frac{|a|^2}{n} \biggr) n^{d/2} \bigl|\pi_n(a)-
p_n(a) \bigr| \biggr) \leq\delta_n.
\end{equation}
Let $(k_n),(\ell_n),(m_n)$ be three sequences of positive integers such
that $n^{-1/2}k_n\rightarrow u$,
$n^{-1/2}\ell_n\rightarrow v$ and $n^{-1/2}m_n\rightarrow w$, where
$0<w<u\wedge v$. Write
\begin{eqnarray*}
&&n^{d/2} \Phi^n_{x_n,y_n}(k_n,
\ell_n,m_n)\\
&&\qquad = n^{3d/4} \int_{\mathbb{R}^d}
\,\mathrm{d}z\, \pi_{m_n}\bigl(\bigl\lfloor zn^{1/4}\bigr\rfloor
\bigr) \pi_{k_n-m_n}\bigl(x_n-\bigl\lfloor zn^{1/4}
\bigr\rfloor\bigr)\pi_{\ell
_n-m_n}\bigl(y_n-\bigl\lfloor
zn^{1/4}\bigr\rfloor\bigr),
\end{eqnarray*}
and note that, for every fixed $z\in\mathbb{R}^d$,
\begin{eqnarray*}
\lim_{n\to\infty} n^{d/4} \pi_{m_n}\bigl(\bigl
\lfloor zn^{1/4}\bigr\rfloor\bigr)&=& p_w(z),
\\
\lim_{n\to\infty} n^{d/4} \pi_{k_n-m_n}
\bigl(x_n-\bigl\lfloor zn^{1/4}\bigr\rfloor\bigr)&=&
p_{u-w}(x-z),
\\
\lim_{n\to\infty} n^{d/4} \pi_{\ell_n-m_n}
\bigl(y_n-\bigl\lfloor zn^{1/4}\bigr\rfloor \bigr)&=&
p_{v-w}(y-z),
\end{eqnarray*}
by \eqref{LLT}. These convergences even hold uniformly in $z$. It then
follows that
%
%
%e12 #&#
\begin{eqnarray}
\label{Tlocal2} \lim_{n\to\infty} n^{d/2}
\Phi^n_{x_n,y_n}(k_n,\ell_n,m_n)
&=& \int_{\mathbb{R}^d}\, \mathrm{d}z\, p_w(z)p_{u-w}(x-z)p_{v-w}(y-z)
\nonumber
\\[-8pt]
\\[-8pt]
\nonumber
&=:& \Psi_{x,y}(u,v,w).
\end{eqnarray}
Indeed, using \eqref{LLT} again, we have, for every
$K>2(|x|\vee|y|)+2$ and every sufficiently large $n$,
\begin{eqnarray*}
&&n^{3d/4} \int_{\{|z|\geq K+1\}} \,\mathrm{d}z\, \pi_{m_n}
\bigl(\bigl\lfloor zn^{1/4}\bigr\rfloor\bigr) \pi_{k_n-m_n}
\bigl(x_n-\bigl\lfloor zn^{1/4}\bigr\rfloor\bigr)
\pi_{\ell
_n-m_n}\bigl(y_n-\bigl\lfloor zn^{1/4}\bigr
\rfloor\bigr)
\\
&&\qquad \leq C \int_{\{|z|\geq K+1\}} \,\mathrm{d}z \biggl(\frac
{1}{(|z|-1)^2}
\biggr)^3,
\end{eqnarray*}
with a constant $C$ independent of $n$ and $K$. The right-hand side of
the last display
tends to $0$ as $K$ tends to infinity. Together with the previously
mentioned uniform
convergence, this suffices to justify \eqref{Tlocal2}.

By \cite{LG1}, Theorem~1.15, we have
\[
\biggl(\frac{\rho}{2} n^{-1/2} H^n_{\lfloor nt\rfloor}
\biggr)_{0\leq
t\leq
1} \build{\la}_{n\to\infty}^{\mathrm{(d)}} (
\mathbf{e}_t )_{0\leq t\leq1},
\]
where $(\mathbf{e}_t)_{0\leq t\leq1}$
is a normalized Brownian excursion, and we recall that $\rho^2$ is the
variance of $\mu$. The latter convergence holds in the sense of the
weak convergence of laws on the Skorokhod space $\mathbb{D}([0,1],\R
_{+})$ of c\`adl\`ag functions from $[0,1]$
into~$\R_+$. Using the Skorokhod representation theorem,
we may and will assume that this convergence holds almost surely, uniformly
in $t\in[0,1]$. Recalling~\eqref{MRCA}, it follows that we have also
\[
\frac{\rho}{2} n^{-1/2} \check H^n_{\lfloor ns\rfloor,\lfloor
nt\rfloor}
\build{\la}_{n\to\infty}^{}\min_{s\wedge t\leq r\leq s\vee t} \mathbf
{e}_r =: m_\mathbf{e}(s,t),
\]
uniformly in $s,t\in[0,1]$, a.s.

As a consequence of \eqref{Tlocal2} and the preceding observations, we have,
for every $s,t\in(0,1)$ with $s\neq t$,
%
%
%e13 #&#
\begin{eqnarray}
\label{Tlocal3} &&\lim_{n\to\infty} n^{d/2}
\Phi^n_{x_n,y_n} \bigl(H^n_{\lfloor
ns\rfloor
},H^n_{\lfloor nt\rfloor},
\check H^n_{\lfloor ns\rfloor,\lfloor nt\rfloor} \bigr)
\nonumber
\\[-8pt]
\\[-8pt]
\nonumber
&&\qquad = \Psi_{x,y} \biggl(
\frac{2}{\rho}\mathbf{e}_s,\frac{2}{\rho
}\mathbf
{e}_t,\frac{2}{\rho}m_\mathbf{e}(s,t) \biggr),\qquad
\mbox{a.s.}
\end{eqnarray}

We claim that we can deduce from \eqref{Tlocal1} and \eqref{Tlocal3} that
%
%
%e14 #&#
\begin{eqnarray}
\label{Tlocal4} &&\lim_{n\to\infty} n^{{d}/{2}-2} E \bigl[
L_n(x_n) L_n(y_n) \bigr]
\nonumber
\\[-8pt]
\\[-8pt]
\nonumber
&&\qquad= E
\biggl[ \int_0^1\int_0^1
\,\mathrm{d}s \,\mathrm{d}t\, \Psi_{x,y} \biggl(\frac{2}{\rho}
\mathbf{e}_s,\frac{2}{\rho}\mathbf{e}_t,
\frac
{2}{\rho
}m_\mathbf{e}(s,t) \biggr) \biggr].
\end{eqnarray}
Note that the right-hand side of \eqref{Tlocal4} coincides with the
function $\varphi(x,y)$
in the lemma. To see this, we can use Theorem III.6 of \cite{Zurich} to
verify that the joint density of the triple
\[
\bigl(m_\mathbf{e}(s,t), \mathbf{e}_s-m_\mathbf{e}(s,t),
\mathbf {e}_t-m_\mathbf{e}(s,t) \bigr)
\]
when $s$ and $t$ are chosen uniformly over $[0,1]$, independently and
independently of $\mathbf{e}$,
is
\[
16(r_1+r_2+r_3)\exp\bigl(-2(r_1+r_2+r_3)^2
\bigr).
\]

So the proof of the second assertion will be complete if we can justify
\eqref{Tlocal4}. By Fatou's lemma,
\eqref{Tlocal1} and \eqref{Tlocal3},
we have first
\[
\liminf_{n\to\infty} n^{{d}/{2}-2} E \bigl[ L_n(x_n)
L_n(y_n) \bigr] \geq E \biggl[ \int_0^1
\int_0^1 \,\mathrm{d}s \,\mathrm{d}t\, \Psi
_{x,y} \biggl(\frac{2}{\rho}\mathbf{e}_s,
\frac{2}{\rho}\mathbf{e}_t,\frac
{2}{\rho
}m_\mathbf{e}(s,t)
\biggr) \biggr].
\]
Furthermore, dominated convergence shows that, for every $K>0$,
\begin{eqnarray*}
&&\lim_{n\to\infty} E \biggl[\int_0^1
\int_0^1 \,\mathrm{d}s \,\mathrm{d}t
\bigl(n^{d/2}\Phi ^n_{x_n,y_n}\bigl(H^n_{\lfloor ns\rfloor},H^n_{\lfloor nt\rfloor},
\check H^n_{\lfloor ns\rfloor,\lfloor nt\rfloor}\bigr)\wedge K \bigr) \biggr]
\\
&&\qquad = E \biggl[\int_0^1\int
_0^1 \,\mathrm{d}s \,\mathrm{d}t \biggl(\Psi
_{x,y} \biggl(\frac{2}{\rho}\mathbf{e}_s,
\frac{2}{\rho}\mathbf {e}_t,\frac
{2}{\rho}m_\mathbf{e}(s,t)
\biggr) \wedge K \biggr) \biggr].
\end{eqnarray*}
Write $\Gamma_n(s,t)= n^{d/2}\Phi^n_{x_n,y_n}(H^n_{\lfloor ns\rfloor
},H^n_{\lfloor nt\rfloor}, \check
H^n_{\lfloor ns\rfloor,\lfloor nt\rfloor})$ to simplify notation. In
view of the preceding comments, it will be enough to verify that
%
%
%e15 #&#
\begin{equation}
\label{Tlocal-key} \lim_{K\to\infty} \biggl( \limsup
_{n\to\infty} E \biggl[\int_0^1\int
_0^1 \,\mathrm{d}s \,\mathrm{d}t\,
\Gamma_n(s,t) \mathbf{1}_{\{\Gamma_n(s,t)>K\}} \biggr] \biggr) =0.
\end{equation}
To this end, we will make use of the bound
%
%
%e16 #&#
\begin{equation}
\label{boundpi} \sup_{k\geq0} \pi_k(x)\leq M
\bigl(|x|^{-d}\wedge1\bigr),
\end{equation}
which holds for every $x\in\mathbb{Z}^d$ with a constant $M$
independent of $x$. This bound can be
obtained easily by combining \eqref{LLT} and Proposition~2.4.6 in
\cite{LL}. Then let $k,\ell,m\geq0$
be integers such that $k\wedge\ell\geq m$, and recall that
\[
\Phi^n_{x_n,y_n}(k,\ell,m)= \sum_{a\in\mathbb{Z}^d}
\pi_m(a) \pi _{k-m}(x_n-a)
\pi_{\ell-m}(y_n-a).
\]
Fix $\ve>0$ such that $|x|\wedge|y|>2\ve$.
Consider first the contribution to the sum in the right-hand side
coming from values of
$a$ such that $|a|\leq\ve n^{1/4}$. For such values of~$a$ (and
assuming that $n$ is large enough), the estimate
\eqref{boundpi} allows us to bound both $\pi_{k-m}(x_n-a)$ and $\pi
_{\ell-m}(y_n-a)$ by $M\ve^{-d}n^{-d/4}$.
On the other hand, if $|a|\geq\ve n^{1/4}$, we can bound $\pi_m(a)$ by
$M\ve^{-d}n^{-d/4}$, whereas
\eqref{LLT} shows that
the sum
\[
\sum_{|a|\geq\ve n^{1/4}} \pi_{k-m}(x_n-a)
\pi_{\ell-m}(y_n-a)
\]
is bounded above by $c_1((k-m)^{-d/2}\wedge(\ell-m)^{-d/2}\wedge1)$
for some constant $c_1$. Summarizing, we get the bound
\begin{eqnarray*}
&&\Phi^n_{x_n,y_n}(k,\ell,m)\\
&&\qquad\leq M^2
\ve^{-2d}n^{-d/2} + c_1M\ve ^{-d}n^{-d/4}
\bigl((k-m)^{-d/2}\wedge(\ell-m)^{-d/2}\wedge1\bigr)
\\
&&\qquad\leq c_{1,\ve}n^{-d/2} + c_{2,\ve}n^{-d/4}
\bigl((k+\ell-2m)^{-d/2} \wedge1\bigr),
\end{eqnarray*}
where $c_{1,\ve}$ and $c_{2,\ve}$ are constants that do not depend on
$n,k,\ell,m$. Then observe that, for every
$s,t\in(0,1)$,
\[
H^n_{\lfloor ns\rfloor}+H^n_{\lfloor nt\rfloor}-2 \check
H^n_{\lfloor ns\rfloor,\lfloor nt\rfloor} = d_n \bigl(u^n_{\lfloor
ns\rfloor},u^n_{\lfloor nt\rfloor}
\bigr),
\]
where $d_n$ denotes the usual graph distance on $\mathcal{T}_n$. From
the preceding bound, we thus get
\[
\Gamma_n(s,t) \leq c_{1,\ve} + c_{2,\ve}n^{d/4}
\bigl(d_n\bigl(u^n_{\lfloor
ns\rfloor},u^n_{\lfloor nt\rfloor}
\bigr)^{-d/2}\wedge1 \bigr).
\]
It follows that, for every $K>0$,
\begin{eqnarray*}
&&\int_0^1\int_0^1
\,\mathrm{d}s \,\mathrm{d}t\, \Gamma_n(s,t) \mathbf {1}_{\{\Gamma_n(s,t)>c_{1,\ve}+ c_{2,\ve}K\}}
\\
&&\qquad \leq\int_0^1\int_0^1
\,\mathrm{d}s \,\mathrm{d}t \bigl(c_{1,\ve} + c_{2,\ve}n^{d/4}
\bigl(d_n\bigl(u^n_{\lfloor ns\rfloor},u^n_{\lfloor nt\rfloor
}
\bigr)^{-d/2}\wedge1\bigr) \bigr) \\
&&\qquad\quad{}\times\mathbf{1}_{ \{ n^{d/4}d_n(u^n_{\lfloor ns\rfloor},u^n_{\lfloor
nt\rfloor})^{-d/2} >K \}}
\\
&&\qquad = n^{-2} \sum_{u,v\in\mathcal{T}_n}
\bigl(c_{1,\ve} + c_{2,\ve
}n^{d/4} \bigl(d_n(u,v)^{-d/2}
\wedge1\bigr) \bigr) \mathbf{1}_{ \{ d_n(u,v)< K^{-2/d} n^{1/2} \}}.
\end{eqnarray*}
By an estimate found in Theorem~1.3 of \cite{DJ}, there exists a
constant $c_0$ that only depends
on $\mu$, such that, for every integer $k\geq1$,
%
%
%e17 #&#
\begin{equation}
\label{estidistance} E \bigl[\#\bigl\{(u,v)\in\mathcal{T}_n\times
\mathcal{T}_n\dvtx d_n(u,v)=k\bigr\} \bigr] \leq
c_0kn.
\end{equation}
It then follows that
\begin{eqnarray*}
&&E \biggl[\int_0^1\int_0^1
\,\mathrm{d}s \,\mathrm{d}t\, \Gamma_n(s,t) \mathbf{1}_{\{\Gamma_n(s,t)>c_{1,\ve}+ c_{2,\ve}K\}}
\biggr]
\\
&&\qquad \leq n^{-1}\bigl(c_{1,\ve} + c_{2,\ve}n^{d/4}
\bigr) + c_0 n^{-1} \sum_{k=1}^{\lfloor K^{-2/d} n^{1/2}\rfloor}
k\bigl(c_{1,\ve} + c_{2,\ve} n^{d/4} k^{-d/2}
\bigr).
\end{eqnarray*}
It is now elementary to verify that
the right-hand side of the preceding display has a limit $g(K)$ when
$n\to\infty$, and that $g(K)$ tends to $0$ as $K\to\infty$ (note that
we use the fact
that $d\leq3$). This completes the proof
of \eqref{Tlocal-key} and of the second assertion of the lemma.

The proof of the first assertion is similar and easier. We first note that
\[
E \bigl[L_n(b)^2 \bigr]=E \biggl[ \sum
_{u,v\in\mathcal{T}_n} \Phi ^n_{b,b}\bigl(|u|,|v|,|u\wedge
v|\bigr) \biggr],
\]
where the function $ \Phi^n_{b,b}$ is defined as above. Then, assuming that
$|b|\geq2\ve n^{1/4}$, the same arguments as in the first part of the
proof give the
bound
\[
\Phi^n_{b,b}\bigl(|u|,|v|,|u\wedge v|\bigr)\leq c_{1,\ve}
n^{-d/2} + c_{2,\ve} n^{-d/4} \bigl(d_n(u,v)^{-d/2}
\wedge1\bigr).
\]
By summing over all choices of $u$ and $v$, it follows that
\begin{eqnarray*}
&&E \bigl[L_n(b)^2 \bigr]\\
&&\qquad\leq c_{1,\ve}n^{2-d/2}
\\
&&\qquad\quad{}+ c_{2,\ve
}n^{-d/4} \biggl( n + E \biggl[\sum
_{u,v\in\mathcal{T}_n, 1\leq d_n(u,v)\leq\sqrt{n}} d_n(u,v)^{-d/2} \biggr] +
n^2\times n^{-d/4} \biggr)
\\
&&\qquad\leq(c_{1,\ve}+2c_{2,\ve})n^{2-d/2}\\
&&\qquad\quad{} +
c_{2,\ve}n^{-d/4} \sum_{k=1}^{\lfloor\sqrt{n}\rfloor}
k^{-d/2} E \bigl[\#\bigl\{(u,v)\in\mathcal{T}_n\dvtx
d_n(u,v)=k\bigr\} \bigr],
\end{eqnarray*}
and the bound stated in the first assertion easily follows from \eqref
{estidistance}.

Let us finally establish the continuity of $\varphi$. We fix $\ve>0$ and
verify that $\varphi$ is continuous on the set $\{|x|\geq2 \ve,|y|\geq
2 \ve\}$. We split
the integral in $\mathrm{d}z$ in two parts:

\begin{longlist}[--]
\item[--] The integral over $|z|\leq\ve$. Write $\varphi_{1,\ve}(x,y)$
for the contribution of this integral. We observe that,
if $|z|\leq\ve$, the function $x\mapsto p_{r_2}(x-z)$ is Lipschitz
uniformly in $z$ and in $r_2$ on the set
$\{|x|\geq2 \ve\}$, and a similar property holds for the function
$y\mapsto p_{r_3}(y-z)$. It follows that
$\varphi_{1,\ve}$ is a Lipschitz function of $(x,y)$ on the set $\{
|x|\geq2 \ve,|y|\geq2 \ve\}$.
\item[--] The integral over $|z|>\ve$. Write $\varphi_{2,\ve}(x,y)$
for the contribution of this integral. Note that
if $(u_n,v_n)_{n\geq1}$ is a sequence in $\mathbb{R}^d\times\mathbb
{R}^d$ such that $|u_n|\wedge|v_n|\geq2\ve$ for every $n$, and
$(u_n,v_n)$ converges to $(x,y)$
as $n\to\infty$, we have, for every fixed $r_1,r_2,r_3>0$,
\begin{eqnarray*}
&&\int_{\{|z|>\ve\}} \,\mathrm{d}z\, p_{r_1}(z)p_{r_2}(u_n-z)p_{r_3}(v_n-z)
\\
&&\qquad\build\la_{n\to
\infty}^{} \int_{\{|z|>\ve\}}
\,\mathrm{d}z\, p_{r_1}(z)p_{r_2}(x-z)p_{r_3}(y-z).
\end{eqnarray*}
We can then use dominated convergence, since there exist constants
$c_\ve$ and $\widetilde c_\ve$
that depend only on $\ve$, such that
\[
\int_{\{|z|>\ve\}} \,\mathrm{d}z\, p_{r_1}(z)p_{r_2}(u_n-z)p_{r_3}(v_n-z)
\leq c_\ve p_{r_2+r_3}(u_n-v_n) \leq
\widetilde c_\ve(r_2+r_3)^{-d/2},
\]
and the right-hand side is integrable for the measure
$(r_1+r_2+r_3)\times\break  e^{-\rho^2(r_1+r_2+r_3)^2/2}\,\mathrm{d}r_1\,\mathrm
{d}r_2\,\mathrm{d}r_3$.
It follows that $\varphi_{2,\ve}$ is also continuous on the set $\{
|x|\geq2 \ve,|y|\geq2 \ve\}$.
\end{longlist}
The preceding considerations complete the proof.
\end{pf*}

In what follows, we use the notation $W^{(1)}=(W^{(1)}_s)_{0\leq s\leq1}$
for a process distributed according to $\N^{(1)}_0$. We recall a result
of Janson and Marckert \cite{JM} that will play an important role below.
As in the proof of Lemma~\ref{contiTlocal}, we let $u^n_0,u^n_1,\ldots,u^n_{n-1}$ be the vertices of $\mathcal{T}_n$
listed in lexicographical order. For every $j\in\{0,1,\ldots,n-1\}$
write $Z^n_j=Z^n(u^n_j)$ for the spatial location of $u^n_j$, and $Z^n_n=0$
by convention. Recalling our assumption \eqref{hypoJM}, we get from
\cite{JM}, Theorem~2, that
%
%
%e18 #&#
\begin{equation}
\label{JaMa} \biggl(\sqrt{\frac{\rho}{2}} n^{-1/4}
Z^n_{\lfloor nt\rfloor} \biggr)_{0\leq t\leq1} \build{
\la}_{n\to\infty}^{\mathrm{(d)}} \bigl(M_\theta ^{1/2}
\widehat W^{(1)}_t \bigr)_{0\leq t\leq1},
\end{equation}
where as usual $M_\theta^{1/2}$ is the unique positive definite
symmetric matrix such that $M_\theta=(M_\theta^{1/2})^2$,
and the convergence holds in distribution in the Skorokhod space
$\mathbb{D}([0,1],\mathbb{R}^{d})$. Note that there are
two minor differences between \cite{JM} and the present setting. First,
\cite{JM} considers one-dimensional labels, whereas our
spatial locations take values in $\Z^d$. However, we can simply project
$Z_n(u)$ on the coordinate axes to get
tightness in the convergence \eqref{JaMa} from the results of \cite
{JM}, and convergence of finite-dimensional marginals
is easy just as in \cite{JM}, Proof of Theorem~1. Second, the
``discrete snake'' of \cite{JM} lists the labels encountered
when exploring the tree $\mathcal{T}_n$ in depth first traversal (or contour
order), whereas we are here enumerating
the vertices in lexicographical order. Nevertheless, the very same
arguments that are used to relate the contour process and the height
function of a random tree (see~\cite{MaMo} or~\cite{LG1}, Section~1.6)
show that asymptotics for the discrete snakes of \cite{JM}
imply similar asymptotics for the labels listed in lexicographical
order of vertices.

In the next theorem, the notation $(l^x,x\in\R^d)$ stands for the
collection of local times of
$W^{(1)}$, which are defined as the continuous density of the
occupation measure of $W^{(1)}$ as in Proposition~\ref{densityZ}.
We define a constant $c>0$ by setting
%
%
%e19 #&#
\begin{equation}
\label{constant-c} c:= \frac{1}{\sigma} \sqrt{\frac{\rho}{2}},
\end{equation}
where $\sigma^2=(\operatorname{ det} M_\theta)^{1/d}$ as previously. We also use
the notation $M_\theta^{-1/2}=(M_\theta^{1/2})^{-1}$.

%
%th4 #&#
\begin{theorem}
\label{convLT}
Let $x^1,\ldots,x^p \in\mathbb{R}^d\setminus\{0\}$, and let
$(x^1_n),\ldots, (x^p_n)$
be sequences in $\mathbb{Z}^d$ such that $\sqrt{\frac{\rho}{2}}
n^{-1/4} M_\theta^{-1/2}x^j_n\la x^j$ as $n\to\infty$, for every
$1\leq j\leq p$. Then
\[
\bigl(n^{{d}/{4}-1}L_n\bigl(x^1_n\bigr),
\ldots, n^{
{d}/{4}-1}L_n\bigl(x^p_n\bigr)
\bigr)\build{\la}_{n\to\infty}^{\mathrm{(d)}} \bigl(c^dl^{x^1},
\ldots,c^dl^{x^p} \bigr),
\]
where the
constant $c$ is given by \eqref{constant-c}.
\end{theorem}

\begin{remarks*}
(i) As mentioned in the \hyperref[sec1]{Introduction}, this result should
be compared with Theorem~1 in \cite{LZ2}, which deals with local
times of branching random walk in $\Z^d$ for $d=2$ or $3$. See also
\cite{BMJ}, Theorem~3.6 and \cite{DJ}, Theorem~1.1, for stronger versions
of the convergence in Theorem~\ref{convLT} when $d=1$.

(ii) It is likely that the result of Lemma~\ref{contiTlocal}
still holds when $x=0$ or $y=0$, and then the condition $x^i\neq0$
in the preceding theorem could
be removed, using also the remark after Lemma~\ref{contiTlocal}.
Proving this reinforcement of Lemma~\ref{contiTlocal} would however
require additional technicalities. Since this extension is not needed
in the proof of our main results, we will not
address this problem here.
\end{remarks*}

\begin{pf*}{Proof of Theorem \ref{convLT}}
To simplify the presentation, we give the details of the proof only in
the isotropic case where
$M_\theta=\sigma^2 \operatorname{ Id}$ (the nonisotropic case is treated in
exactly the same manner at the
cost of a somewhat heavier notation). Our condition on the sequences
$(x^j_n)$ then just says that $c n^{-1/4} x^j_n\la x^j$ as $n\to\infty$.

By the Skorokhod
representation theorem, we may and will assume that the convergence
\eqref{JaMa} holds a.s. To obtain the result of the
theorem, it is then enough to verify that, if $x \in\mathbb
{R}^d\setminus\{0\}$ and $(x_n)$ is a sequence
in $\mathbb{Z}^d$ such that $c n^{-1/4} x_n\la x$ as $n\to\infty$,
we have
%
%
%e20 #&#
\begin{equation}
\label{convLTtech1} n^{{d}/{4}-1}L_n(x_n) \build{
\la}_{n\to\infty}^{\mathrm{(P)}} c^d l^x.
\end{equation}
To this end, fix $x$ and the sequence $(x_n)$, and for every $\ve\in
(0,|x|)$, let $g_\ve$
be a nonnegative continuous function on $\mathbb{R}^d$, with compact
support contained in the open ball of radius $\ve$
centered at $x$, and such that
\[
\int_{\R^{d}} g_\ve(y) \,\mathrm{d}y = 1.
\]
It follows from \eqref{JaMa} (which we assume to hold a.s.) that, for
every fixed $\ve\in(0,|x|)$,
\[
\int_0^1 g_\ve\bigl(c
n^{-1/4} Z^n_{\lfloor nt\rfloor}\bigr) \,\mathrm{d}t \build {
\la}_{n\to\infty}^{\mathrm{ a.s.}} \int_0^1
g_\ve\bigl(\widehat W^{(1)}_t\bigr)\,
\mathrm{d}t.
\]
Furthermore,
\[
\int_0^1 g_\ve\bigl(\widehat
W^{(1)}_t\bigr) \,\mathrm{d}t =\int_{\R^{d}}
g_\ve (y) l^y \,\mathrm{d}y \build{\la}_{\ve\to0}^{\mathrm{ a.s.}}
l^x,
\]
by the continuity of local times.
Let $\delta>0$. By combining the last two convergences, we can find
$\ve
_1\in(0,|x|)$
such that, for every $\ve\in(0,\ve_1)$, there exists an integer
$n_1(\ve
)$ so that for every $n\geq n_1(\ve)$,
%
%
%e21 #&#
\begin{equation}
\label{convLTtech2} P \biggl( \biggl|\int_0^1
g_\ve\bigl(c n^{-1/4} Z^n_{\lfloor nt\rfloor}\bigr)\,
\mathrm{d}t - l^x \biggr| > \delta \biggr) <\delta.
\end{equation}

However, we have
\begin{eqnarray*}
\int_0^1 g_\ve
\bigl(cn^{-1/4}Z^n_{\lfloor nt\rfloor}\bigr)\, \mathrm{d}t &=&
\frac
{1}{n} \sum_{a\in\mathbb{Z}^d} g_\ve
\bigl(cn^{-1/4}a\bigr) L_n(a)
\\
&=& n^{{d}/{4} -1} \int_{\mathbb{R}^d} g_\ve
\bigl(cn^{-1/4}\bigl\lfloor n^{1/4}y\bigr\rfloor\bigr)
L_n\bigl(\bigl\lfloor n^{1/4}y\bigr\rfloor\bigr)\,
\mathrm{d}y.
\end{eqnarray*}
Set
\[
\eta_n(\ve):= \int_{\mathbb{R}^d} g_\ve
\bigl(cn^{-1/4}\bigl\lfloor n^{1/4}y\bigr\rfloor\bigr)
\,\mathrm{d}y
\]
and note that
\[
\eta_n(\ve)\build{\la}_{n\to\infty}^{} \int
_{\mathbb{R}^d} g_\ve(cy) \,\mathrm{d}y = c^{-d}.
\]
By the Cauchy--Schwarz inequality,
\begin{eqnarray*}
&&E \biggl[ \biggl(\int_0^1 g_\ve
\bigl(cn^{-1/4}Z^n_{\lfloor nt\rfloor}\bigr) \,\mathrm {d}t -
\eta_n(\ve) n^{{d}/{4}-1} L_n(x_n)
\biggr)^2 \biggr]
\\
&&\qquad =E \biggl[ \biggl(n^{{d}/{4}-1} \int_{\mathbb{R}^d}
g_\ve \bigl(cn^{-1/4}\bigl\lfloor n^{1/4}y\bigr
\rfloor\bigr) \bigl(L_n\bigl(\bigl\lfloor n^{1/4}y\bigr
\rfloor\bigr) -L_n(x_n)\bigr)\, \mathrm{d}y
\biggr)^2 \biggr]
\\
& &\qquad\leq\eta_n(\ve)\times n^{{d}/{2}-2} \int_{\mathbb{R}^d}
\,\mathrm{d}y\, g_\ve\bigl(cn^{-1/4}\bigl\lfloor
n^{1/4}y\bigr\rfloor\bigr) E \bigl[\bigl(L_n\bigl(\bigl
\lfloor n^{1/4}y\bigr\rfloor\bigr) -L_n(x_n)
\bigr)^2 \bigr].
\end{eqnarray*}
Using the first assertion of Lemma~\ref{contiTlocal}, one easily gets
that, for every fixed $\ve\in(0,|x|)$,
\[
n^{{d}/{2}-2} \int_{\mathbb{R}^d} \,\mathrm{d}y \bigl|g_\ve
\bigl(cn^{-1/4}\bigl\lfloor n^{1/4}y\bigr\rfloor
\bigr)-g_\ve(cy)\bigr | E \bigl[\bigl(L_n\bigl(\bigl\lfloor
n^{1/4}y\bigr\rfloor\bigr) -L_n(x_n)
\bigr)^2 \bigr] \build{\la}_{n\to\infty}^{} 0.
\]
On the other hand, by the second assertion of the lemma,
\begin{eqnarray*}
&&n^{{d}/{2}-2} \int_{\mathbb{R}^d} \,\mathrm{d}y\, g_\ve(cy)
E \bigl[\bigl(L_n\bigl(\bigl\lfloor n^{1/4}y\bigr\rfloor
\bigr) -L_n(x_n)\bigr)^2 \bigr]\\
&&\qquad \build{
\la}_{n\to\infty}^{} \int_{\mathbb{R}^d}\, \mathrm{d}y\,
g_\ve(cy) \biggl(\varphi(y,y)-2\varphi\biggl(\frac
{x}{c},y
\biggr)+\varphi\biggl(\frac{x}{c},\frac{x}{c}\biggr) \biggr).
\end{eqnarray*}
If $\gamma_\ve$ stands for the limit in the last display, the
continuity of $\varphi$
ensures that $\gamma_\ve$ tends to $0$ as $\ve\to0$.

From the preceding considerations, we have
\[
\limsup_{n\to\infty} E \biggl[ \biggl( \int_0^1
g_\ve \bigl(cn^{-1/4}Z^n_{\lfloor
nt\rfloor}\bigr)\,
\mathrm{d}t - \eta_n(\ve) n^{{d}/{4}-1} L_n(x_n)
\biggr)^2 \biggr] \leq c^{-d} \gamma_\ve.
\]
Hence, we can find $\ve_2\in(0,|x|)$ small enough so that, for every
$\ve\in(0,\ve_2)$, there exists an integer
$n_2(\ve)$ such that, for every $n\geq n_2(\ve)$,
%
%
%e22 #&#
\begin{equation}
\label{convLTtech3} P \biggl(\biggl | \int_0^1
g_\ve\bigl(cn^{-1/4}Z^n_{\lfloor nt\rfloor}\bigr)\,
\mathrm {d}t - \eta_n(\ve) n^{{d}/{4}-1} L_n(x_n)
\biggr| > \delta \biggr) <\delta.
\end{equation}
By combining \eqref{convLTtech2} and \eqref{convLTtech3}, we see that,
for every $\ve\in(0,\ve_1\wedge\ve_2)$
and $n\geq n_1(\ve)\vee n_2(\ve)$,
\[
P \bigl( \bigl|\eta_n(\ve) n^{{d}/{4}-1} L_n(x_n)
- l^x \bigr| > 2\delta \bigr) < 2\delta.
\]
Our claim \eqref{convLTtech1} easily follows, since $\eta_n(\ve)$ tends
to $c^{-d}$ as $n\to\infty$.
\end{pf*}

Set $R_n=\#\{ Z^n(u) \dvtx  u\in\mathcal{T}_n\}$. Recall the constant
$c$ from \eqref{constant-c}, and also recall that $\lambda_d$
denotes Lebesgue measure on $\R^d$.

%
%th5 #&#
\begin{theorem}
\label{conv-range}
We have
\[
n^{-d/4} R_n \build{\la}_{n\to\infty}^{\mathrm{(d)}}
c^{-d}\lambda _d(\mathcal{S}),
\]
where $\mathcal{S}$ stands for the support of ISE.
\end{theorem}

\begin{pf}
Again, for the sake of simplicity, we give details only in the
isotropic case $M_\theta=\sigma^2 {\mathrm {Id}}$.
From the definition of ISE, we may take $\mathcal{S}=  \{\widehat
W^{(1)}_t \dvtx 0\leq t\leq1 \}$. We then set, for every $\ve>0$,
\[
\mathcal{S}_\ve:= \bigl\{x\in\R^d \dvtx  \operatorname{dist}(x,
\mathcal{S})\leq \ve \bigr\}.
\]
As in the preceding proof, we may and will assume that the convergence
\eqref{JaMa} holds
almost surely. It then follows that, for every $\ve>0$,
\[
P \bigl(\bigl\{ cn^{-1/4} Z^n(u) \dvtx  u\in
\mathcal{T}_n\bigr\}\subset\mathcal {S}_\ve \bigr) \build{
\la}_{n\to\infty}^{} 1.
\]
Fix $K>0$, and let $B(0,K)$ stand for the closed ball of radius $K$
centered at $0$ in~$\R^d$. Also set
$\mathcal{S}_\ve^{(K)}:= \mathcal{S}_\ve\cap B(0,K+\ve)$. It follows
that we have also
\[
P \bigl(\bigl(\bigl\{ Z^n(u) \dvtx  u\in\mathcal{T}_n
\bigr\} \cap B\bigl(0, c^{-1}n^{1/4} K\bigr)\bigr)\subset
c^{-1}n^{1/4} \mathcal{S}_\ve ^{(K)}
\bigr) \build{\la}_{n\to\infty}^{} 1.
\]
Applying the latter convergence with $\ve$ replaced by $\ve/2$, we get
\[
P \bigl(\#\bigl(\bigl\{ Z^n(u) \dvtx  u\in\mathcal{T}_n
\bigr\} \cap B\bigl(0, c^{-1}n^{1/4} K\bigr)\bigr) \leq
c^{-d}n^{d/4} \lambda_d\bigl(\mathcal
{S}_\ve^{(K)}\bigr) \bigr)\build{\la}_{n\to\infty}^{}
1.
\]
Write $R_n^{(K)}:= \#(\{ Z^n(u) \dvtx  u\in\mathcal{T}_n\}
\cap B(0, c^{-1}n^{1/4} K))$. Since $\lambda_d(\mathcal{S}_\ve^{(K)})
\downarrow\lambda_d(\mathcal{S}\cap B(0,K))$ as
$\ve\downarrow0$, we obtain that, for every $\delta>0$,
\[
P \bigl(n^{-d/4} R_n^{(K)}\leq c^{-d}
\lambda_d\bigl(\mathcal{S}\cap B(0,K)\bigr) +\delta \bigr)\build{
\la}_{n\to\infty}^{} 1,
\]
and, therefore, since the variables $n^{-d/4}R_n^{(K)}$ are uniformly bounded,
%
%
%e23 #&#
\begin{equation}
\label{conv-rangetech1} \lim_{n\to\infty} E \bigl[ \bigl(n^{-d/4}
R_n^{(K)}- c^{-d}\lambda _d\bigl(
\mathcal{S}\cap B(0,K)\bigr) \bigr)^+ \bigr]= 0.
\end{equation}

On the other hand, we claim that we have also
%
%
%e24 #&#
\begin{equation}
\label{conv-rangetech2} \liminf_{n\to\infty} E \bigl[n^{-d/4}
R_n^{(K)} \bigr] \geq c^{-d}E \bigl[
\lambda_d\bigl(\mathcal{S}\cap B(0,K)\bigr) \bigr].
\end{equation}
To see this, observe that
\begin{eqnarray*}
E \bigl[R_n^{(K)} \bigr]& =& \sum
_{a\in\Z^d \cap B(0, c^{-1}n^{1/4} K)} P\bigl(L_n(a)>0\bigr)
\\
& =& \int_{B(0, c^{-1}n^{1/4} K)} \,\mathrm{d}x\, P\bigl(L_n\bigl(
\lfloor x\rfloor \bigr)>0\bigr) + O\bigl(n^{(d-1)/4}\bigr)
\\
& =& n^{d/4} \int_{B(0, c^{-1} K)} \,\mathrm{d}y\, P
\bigl(L_n\bigl(\bigl\lfloor n^{1/4}y\bigr\rfloor\bigr)>0
\bigr) + O\bigl(n^{(d-1)/4}\bigr)
\end{eqnarray*}
as $n\to\infty$. By Theorem~\ref{convLT}, for every $y\neq0$,
\[
\liminf_{n\to\infty} P\bigl(L_n\bigl(\bigl\lfloor
n^{1/4}y\bigr\rfloor\bigr)>0\bigr) \geq P\bigl(l^{cy}>0\bigr)=
P(cy\in\mathcal{S}),
\]
where the equality is derived from Proposition~\ref{hitting-point}.
Fatou's lemma then gives
\[
\liminf_{n\to\infty} n^{-d/4}E \bigl[R_n^{(K)}
\bigr] \geq\int_{B(0,
c^{-1} K)} \,\mathrm{d}y\, P(cy\in\mathcal{S}) =
c^{-d} E \bigl[\lambda_d\bigl(\mathcal{S}\cap B(0,K)\bigr)
\bigr],
\]
which completes the proof of \eqref{conv-rangetech2}.

Using the trivial identity $|x|=2x^+-x$ for every real $x$,
we deduce from \eqref{conv-rangetech1} and \eqref{conv-rangetech2} that
\[
\lim_{n\to\infty} E \bigl[\bigl |n^{-d/4} R_n^{(K)}-
c^{-d}\lambda _d\bigl(\mathcal{S}\cap B(0,K)\bigr) \bigr|
\bigr]=0.
\]
However, we see from \eqref{JaMa} that, for every $\delta>0$, we can
choose $K$
sufficiently large so that we have both $P(\mathcal{S}\subset
B(0,K))\geq1-\delta$
and $P(R_n^{(K)}=R_n)\geq1-\delta$ for every integer $n$. It then follows
from the previous convergence that $n^{-d/4} R_n$
converges in probability to $c^{-d}\lambda_d(\mathcal{S})$ as $n\to
\infty$, and this completes the proof
of Theorem~\ref{conv-range}.
\end{pf}

%s4 #&#
\section{Branching random walk}\label{sec4}

We will now discuss similar results for branching random walk in
$\Z^d$. We consider a system of particles in $\Z^d$ that evolves
in discrete time
in the following way. At time $n=0$, there are $p$ particles all located
at the origin of $\Z^d$ (we will comment on more general
initial configurations in Section~\ref{bps}). A particle located at
the site
$a\in\Z^d$ at time $n$ gives rise at time $n+1$ to a random number of
offspring
distributed according to $\mu$, and their locations are obtained
by adding to $a$ (independently for each offspring) a spatial displacement
distributed according to $\theta$.

In a more formal way, we consider $p$ independent random spatial trees
\[
\bigl(\mathcal{T}^{(1)},\bigl(Z^{(1)}(u)\bigr)_{u\in\mathcal{T}^{(1)}}
\bigr),\ldots, \bigl(\mathcal{T} ^{(p)},\bigl(Z^{(p)}(u)
\bigr)_{u\in\mathcal{T}^{(p)}} \bigr)
\]
distributed according to $\Pi^*_{\mu,\theta}$, and, for every integer
$n\geq0$,
we consider the random point measure
\[
X^{[p]}_n:= \sum_{j=1}^p
\biggl(\sum_{u\in\mathcal{T}^{(j)},|u|=n} \delta _{Z^{(j)}(u)} \biggr),
\]
which corresponds to the sum of the Dirac point masses at the positions of
all particles alive at time $n$.

The set $\mathcal{V}^{[p]}$ of all sites visited by the particles is
the union over all $n\geq0$ of
the supports of $X^{[p]}_n$, or equivalently
\[
\mathcal{V}^{[p]}= \bigl\{a\in\Z^d \dvtx
a=Z^{(j)}(u)\mbox{ for some }j\in\{1,\ldots,p\}\mbox{ and }u\in
\mathcal{T}^{(j)} \bigr\}.
\]
In a way similar to Theorem~\ref{conv-range}, we are interested in limit
theorems for $\#\mathcal{V}^{[p]}$ when $p\to\infty$. To this end, we
will first state
an analog of the convergence \eqref{JaMa}. For every $j\in\{1,\ldots,p\}
$, let
\[
\varnothing=u^{(j)}_0\prec u^{(j)}_1
\prec\cdots\prec u^{(j)}_{\#
\mathcal{T}^{(j)}-1}
\]
be the vertices
of $\mathcal{T}^{(j)}$ listed in lexicographical order, and set
$H^{(j)}_i=|u^{(j)}_i|$
and $Z^{(j)}_i=Z^{(j)}(u^{(j)}_i)$, for $0\leq i\leq\#\mathcal{T}^{(j)}-1$.
Define the height function
$(H^{[p]}_k,k\geq0)$
of $\mathcal{T}^{(1)},\ldots,\mathcal{T}^{(p)}$ by concatenating the
discrete functions
$(H^{(j)}_i,0\leq i\leq\#\mathcal{T}^{(j)}-1$),
and setting $H^{[p]}_k=0$ for $k\geq\#\mathcal{T}^{(1)}+\cdots+\#
\mathcal{T}^{(p)}$.
Similarly, define the function
$(Z^{[p]}_k,k\geq0)$
by concatenating the discrete functions $(Z^{(j)}_i,0\leq i\leq\#
\mathcal{T}^{(j)}-1$),
and setting $Z^{[p]}_k=0$ for $k\geq\#\mathcal{T}^{(1)}+\cdots+\#
\mathcal{T}^{(p)}$.
Finally, we use linear interpolation to define $H^{[p]}_t$ and
$Z^{[p]}_t$ for every
real $t\geq0$. We can now state our analog of~\eqref{JaMa}.

%
%pr6 #&#
\begin{proposition}
\label{JaMabis}
We have
\begin{eqnarray*}
&&\biggl( \biggl(\frac{\rho}{2} p^{-1} H^{[p]}_{p^2s},
\sqrt{\frac{\rho
}{2}} p^{-{1}/{2}} Z^{[p]}_{p^2s}
\biggr)_{ s\geq0}, p^{-2} \bigl(\#\mathcal{T}^{(1)}+
\cdots+\#\mathcal{T}^{(p)}\bigr) \biggr) \\
&&\qquad\build{\la}_{p\to\infty}^{\mathrm{(d)}}
\bigl(\bigl(\zeta_{s\wedge\tau}, M_\theta ^{1/2}\wh
W_{s\wedge\tau}\bigr)_{s\geq0},\tau \bigr),
\end{eqnarray*}
where $(W_s)_{s\geq0}$ is a standard Brownian snake, $\tau$ denotes
the first
hitting time of $2/\rho$ by the local time at $0$ of the lifetime
process $(\zeta_s)_{s\geq0}$, and the convergence
of processes holds in the sense of the topology of uniform convergence
on compact sets.
\end{proposition}

The joint convergence of the processes $\frac{\rho}{2} p^{-1}
H^{[p]}_{p^2s}$ and of the random variables $p^{-2}
(\#\mathcal{T}^{(1)}+\cdots+\#\mathcal{T}^{(p)})$
is a consequence of \cite{LG1}, Theorem~1.8, see in
particular~(7) and (9) in \cite{LG1} (note that the local times of the
process $(\zeta_s)_{s\geq0}$ are chosen to be right-continuous
in the space variable, so that our local time at $0$ is twice the local
time that appears in the display (7) in \cite{LG1}). Given the latter
joint convergence, the desired statement
can be obtained by following the arguments of the proof of Theorem~2 in
\cite{JM}. The fact that
we are dealing with unconditioned trees makes things easier than in
\cite{JM} and we omit the details.

We now state an intermediate result, which is of independent interest. Under
the probability measure $\Pi^*_{\mu,\theta}$, we let $\bR:= \{z_u
\dvtx  u\in\mathcal{T}\}$
be the set of all points visited by the tree-indexed random walk.

%
%th7 #&#
\begin{theorem}
\label{estim-visit}
We have
\[
\lim_{|a|\to\infty}\bigl |M_\theta^{-1/2}a\bigr|^2
\Pi^*_{\mu,\theta}(a\in \bR) = \frac{2(4-d)}{\rho^2}.
\]
\end{theorem}

\begin{pf} We start by proving the upper bound
\[
\limsup_{|a|\to\infty} \bigl|M_\theta^{-1/2}a\bigr|^2
\Pi^*_{\mu,\theta
}(a\in\bR ) \leq\frac{2(4-d)}{\rho^2}.
\]
By an easy compactness argument, it is enough to prove that, if $(a_k)$
is a sequence
in $\Z^d$ such that $|a_k|\to\infty$ and $a_k/|a_k| \to x$, with
$x\in\R
^d$ and $|x|=1$, then
%
%
%e25 #&#
\begin{equation}
\label{estim-v0} \limsup_{k\to\infty} |a_k|^2
\Pi^*_{\mu,\theta}(a_k\in\bR) \leq \frac
{2(4-d)}{\rho^2|M_\theta^{-1/2}x|^2}.
\end{equation}
Set $p_k=|a_k|^2\in\Z_+$ to simplify notation. We note that
%
%
%e26 #&#
\begin{equation}
\label{estim-v1} P \bigl(a_k\in\mathcal{V}^{[p_k]} \bigr) = 1
- \bigl( 1- \Pi^*_{\mu,\theta
}(a_k\in\bR) \bigr)^{p_k}.
\end{equation}
On the other hand, it follows from our definitions that
\[
P \bigl(a_k\in\mathcal{V}^{[p_k]} \bigr) \leq P \biggl(
\exists s\geq0 \dvtx  \frac{1}{\sqrt{p_k}} Z^{(p_k)}_s =
\frac{a_k}{|a_k|} \biggr).
\]
We can then use Proposition~\ref{JaMabis} to get
\begin{eqnarray*}
\limsup_{k\to\infty} P \bigl(a_k\in
\mathcal{V}^{[p_k]} \bigr) &\leq&\P_0 \biggl(\exists s\in[0,
\tau]\dvtx  M_\theta^{1/2}\wh W_s= \sqrt{
\frac{\rho}{2}} x \biggr)
\\
&=& 1-\exp \biggl(-\frac{2}{\rho} \N_0 \biggl(\sqrt{
\frac{\rho
}{2}}M_\theta ^{-1/2}x\in\mathcal{R} \biggr)
\biggr)
\\
&=& 1-\exp \biggl(- \frac{2(4-d)}{\rho^2|M_\theta^{-1/2}x|^2} \biggr).
\end{eqnarray*}
The second line follows from excursion theory for the Brownian snake,
and the third one
uses the formula for $\N_0(y\in\mathcal{R})$, which has been recalled
already in the proof of
Proposition~\ref{hitting-point}. By combining the bound of the last
display with \eqref{estim-v1},
we get our claim \eqref{estim-v0}, and this completes the proof of the
upper bound.

Let us turn to the proof of the lower bound. As in the proof of the
upper bound, it is enough
to consider a sequence $(a_k)$
in $\Z^d$ such that $|a_k|\to\infty$ and $a_k/|a_k| \to x$, with
$x\in\R
^d$ and $|x|=1$, and then
to verify that
%
%
%e27 #&#
\begin{equation}
\label{estim-v3} \liminf_{k\to\infty} |a_k|^2
\Pi^*_{\mu,\theta}(a_k\in\bR) \geq \frac
{2(4-d)}{\rho^2|M_\theta^{-1/2}x|^2}.
\end{equation}
As previously, we set $p_k=|a_k|^2$. We fix $0<\ve<M$, and we introduce
the function $g_\mu$ defined on $\Z_+$
by $g_\mu(j)=\Pi_\mu(\#\mathcal{T}=j)$. Then
\begin{eqnarray*}
|a_k|^2 \Pi^*_{\mu,\theta}(a_k\in\bR) &
\geq& p_k^3 \int_\ve^M\,
\mathrm{d}r\, \Pi^*_{\mu,\theta}\bigl(a_k\in\bR, \# \mathcal{T} =
\bigl\lfloor p_k^2r\bigr\rfloor\bigr)
\\
&=& p_k^3 \int_\ve^M\,
\mathrm{d}r\, g_\mu \bigl(\bigl\lfloor p_k^2r
\bigr\rfloor \bigr) P \bigl( L_{\lfloor p_k^2 r\rfloor}(a_k)>0 \bigr),
\end{eqnarray*}
where we use the same notation as in Lemma~\ref{contiTlocal}: $L_n(b)$
denotes the number of visits of site $b$ by a random walk indexed by
a tree distributed according to $\Pi_\mu(\cdot\mid\#\mathcal
{T}=n)$. Note
that Theorem~\ref{convLT}
gives, for every $r\in[\ve,M]$,
\[
\liminf_{k\to\infty} P \bigl( L_{\lfloor p_k^2 r\rfloor}(a_k)>0
\bigr) \geq P \bigl(l^{r^{-1/4}z}>0 \bigr),
\]
where we write $z=\sqrt{\frac{\rho}{2}} M_\theta^{-1/2}x$ to
simplify notation.
To complete the argument, we consider for simplicity the aperiodic case
where $\mu$
is not supported on a strict subgroup of $\Z$ [the reader will easily
be able to
extend our method to the general case, using \eqref{Kemp1} instead of
\eqref{Kemp2}]. By
\eqref{Kemp2}, we have for every $r\in[\ve,M]$,
\[
\lim_{k\to\infty} p_k^3 g_\mu
\bigl(\bigl\lfloor p_k^2r\bigr\rfloor \bigr)=
\frac
{1}{\rho\sqrt{2\pi r^3}}.
\]
Using this together with the preceding display, and applying Fatou's
lemma, we obtain
%
%
%e28 #&#
\begin{equation}
\label{estim-v4} \liminf_{k\to\infty} |a_k|^2
\Pi^*_{\mu,\theta}(a_k\in\bR) \geq\int_\ve^M
\frac{\mathrm{d}r}{\rho\sqrt{2\pi r^3}} P \bigl(l^{r^{-1/4}z}>0 \bigr).
\end{equation}
A scaling argument
shows that
\[
P \bigl(l^{r^{-1/4}z}>0 \bigr)= \N^{(1)}_0 \bigl(
\ell^{r^{-1/4}z}>0 \bigr)= \N ^{(r)}_0 \bigl(
\ell^{z}>0 \bigr).
\]
Using this remark and formula \eqref{decoIto}, we see that the
right-hand side of \eqref{estim-v4} can be rewritten as $\frac
{2}{\rho}
\N_0(\mathbf{1}_{\{\ve<\gamma<M\}} \mathbf{1}_{\{\ell^{z}>0\}})$.
By choosing $\ve$ small enough and $M$ large enough, the latter
quantity can be made
arbitrarily close to
\[
\frac{2}{\rho} \N_0\bigl(\ell^{z}>0\bigr)=
\frac{2}{\rho} \biggl(2-\frac{d}{2}\biggr) |z|^{-2}=
\frac{2(4-d)}{\rho^2|M_\theta^{-1/2}x|^2}.
\]
This completes the proof of the lower bound and of Theorem~\ref{estim-visit}.
\end{pf}

Recall our notation $\mathcal{V}^{[p]}$ for the set of all sites
visited by
the branching random walk starting with $p$ initial particles located
at the origin.

%
%th8 #&#
\begin{theorem}
\label{rangeBRW}
We have
\[
p^{-d/2} \#\mathcal{V}^{[p]} \build{\la}_{p\to\infty}^{\mathrm{(d)}}
\biggl(\frac{2\sigma}{\rho}\biggr)^{d} \lambda_d \biggl(
\bigcup_{t\geq0} \operatorname{ supp} X_t \biggr),
\]
where $(X_t)_{t\geq0}$ is a $d$-dimensional super-Brownian motion with
branching mechanism $\psi(u)=2u^2$
started from $\delta_0$, and $\operatorname{ supp} X_t$ denotes the topological
support of~$X_t$.
\end{theorem}

\begin{pf}
Via the Skorokhod representation theorem,
we may and will assume that the convergence in Proposition~\ref
{JaMabis} holds a.s., and
we will then prove that the convergence of the theorem holds in probability.
If $\ve>0$ is fixed, the (a.s.) convergence in Proposition~\ref
{JaMabis} implies that, a.s. for all large enough $p$, we have
\[
\sqrt{\frac{\rho}{2}} p^{-1/2}\mathcal{V}^{[p]} \subset
\mathcal {U}_\ve \bigl(\bigl\{M_\theta^{1/2}\wh
W_s\dvtx 0\leq s\leq\tau\bigr\} \bigr),
\]
where, for any compact subset $\mathcal{K}$ of $\R^d$, $\mathcal
{U}_\ve
(\mathcal{K})$ denotes the set of all
points whose distance from $\mathcal{K}$ is strictly less than $\ve$.
It follows that we have a.s.
\[
\limsup_{p\to\infty} p^{-d/2}\#\mathcal{V}^{[p]}
\leq \biggl(\frac
{2}{\rho} \biggr)^{d/2} \lambda_d \bigl(
\mathcal{U}_{2\ve} \bigl(\bigl\{ M_\theta ^{1/2}\wh
W_s \dvtx 0\leq s\leq\tau\bigr\} \bigr) \bigr).
\]
Since $\ve$ was arbitrary, we also get a.s.
%
%
%e29 #&#
\begin{equation}
\label{BRW1} \limsup_{p\to\infty} p^{-d/2}\#
\mathcal{V}^{[p]}\leq \biggl(\frac
{2}{\rho} \biggr)^{d/2}
\lambda_d \bigl( \bigl\{M_\theta^{1/2}\wh
W_s \dvtx  0\leq s\leq\tau \bigr\} \bigr).
\end{equation}

To get an estimate in the reverse direction, we argue in a way very
similar to the proof of Theorem~\ref{conv-range}.
We fix $K>0$, and note that a minor modification of the preceding
arguments also gives a.s.
\begin{eqnarray*}
&&\limsup_{p\to\infty} p^{-d/2}\# \bigl(\mathcal{V}^{[p]}
\cap B\bigl(0,p^{1/2}K\bigr) \bigr)\\
&&\qquad\leq \biggl(\frac{2}{\rho}
\biggr)^{d/2} \lambda _d \bigl( \bigl\{M_\theta^{1/2}
\wh W_s \dvtx 0\leq s\leq\tau \bigr\}\cap B\bigl(0,K'
\bigr) \bigr),
\end{eqnarray*}
where $K'=\sqrt{\frac{\rho}{2}} K$.
Since the variables $p^{-d/2} \#(\mathcal{V}^{[p]}\cap B(0,p^{1/2}K))$
are uniformly bounded, it follows that
%
%
%e30 #&#
\begin{eqnarray}
\label{BRW2}&& \lim_{p\to\infty} E \biggl[ \biggl(p^{-d/2}
\# \bigl(\mathcal {V}^{[p]}\cap B\bigl(0,p^{1/2}K\bigr) \bigr)\nonumber\\
&&\hspace*{25pt}\qquad{}-
\biggl(\frac{2}{\rho} \biggr)^{d/2} \lambda_d \bigl(\bigl
\{M_\theta^{1/2}\wh W_s\dvtx 0\leq s\leq\tau\bigr\}
\cap B\bigl(0,K'\bigr) \bigr) \biggr)^+ \biggr]\\
&&\qquad=0.\nonumber
\end{eqnarray}
On the other hand,
\begin{eqnarray*}
&&p^{-d/2} E \bigl[\#\bigl(\mathcal{V}^{[p]}\cap B
\bigl(0,p^{1/2}K\bigr)\bigr) \bigr]\\
&&\qquad= p^{-d/2} \sum
_{a\in\Z^d\cap B(0,p^{1/2}K)} P\bigl(a\in\mathcal {V}^{[p]}\bigr)
\\
&&\qquad=p^{-d/2} \sum_{a\in\Z^d\cap B(0,p^{1/2}K)} \bigl(1-\bigl(1-
\Pi^*_{\mu,\theta}(a\in\bR)\bigr)^p \bigr)
\\
&&\qquad\build{\la}_{p\to\infty}^{} \int_{B(0,K)}\,
\mathrm{d}x \biggl(1-\exp \biggl(-\frac{2(4-d)}{\rho^2 |M_\theta^{-1/2}x|^2} \biggr) \biggr),
\end{eqnarray*}
where the last line is an easy consequence of Theorem~\ref{estim-visit}.
Furthermore,
\begin{eqnarray*}
&&E \biggl[ \biggl(\frac{2}{\rho} \biggr)^{d/2}
\lambda_d \bigl(\bigl\{M_\theta ^{1/2}\wh
W_s\dvtx 0\leq s\leq\tau\bigr\}\cap B\bigl(0,K'\bigr)
\bigr) \biggr]
\\
&&\qquad = \biggl(\frac{2}{\rho} \biggr)^{d/2}\int_{B(0,K')}\,
\mathrm{d}y \biggl(1-\exp \biggl(-\frac{2}{\rho} \N_0
\bigl(M_\theta^{-1/2}y\in\mathcal{R}\bigr) \biggr) \biggr)
\\
&&\qquad = \biggl(\frac{2}{\rho} \biggr)^{d/2}\int_{B(0,K')}\,
\mathrm{d}y \biggl(1-\exp \biggl(-\frac{4-d}{\rho|M_\theta^{-1/2}y|^2} \biggr) \biggr)
\\
&&\qquad =\int_{B(0,K)} \,\mathrm{d}x \biggl(1-\exp \biggl(-
\frac{2(4-d)}{\rho
^2 |M_\theta^{-1/2}x|^2} \biggr) \biggr).
\end{eqnarray*}
From the last two displays, we get
%
%
%e31 #&#
\begin{eqnarray}
\label{BRW3} &&\lim_{p\to\infty} E \bigl[p^{-d/2} \#\bigl(
\mathcal{V}^{[p]}\cap B\bigl(0,p^{1/2}K\bigr)\bigr) \bigr]
\nonumber
\\[-8pt]
\\[-8pt]
\nonumber
&&\qquad = E
\biggl[ \biggl(\frac{2}{\rho} \biggr)^{d/2} \lambda_d
\bigl(\bigl\{M_\theta ^{1/2}\wh W_s \dvtx 0\leq s
\leq\tau\bigr\}\cap B\bigl(0,K'\bigr) \bigr) \biggr].
\end{eqnarray}
From \eqref{BRW2} and \eqref{BRW3}, we have
\begin{eqnarray*}
&&\lim_{p\to\infty} E \biggl[\biggl |p^{-d/2} \#\bigl(
\mathcal{V}^{[p]}\cap B\bigl(0,p^{1/2}K\bigr)\bigr)\\
&&\hspace*{20pt}\qquad{}- \biggl(
\frac{2}{\rho} \biggr)^{d/2} \lambda_d \bigl(\bigl
\{M_\theta^{1/2}\wh W_s\dvtx 0\leq s\leq\tau\bigr\}
\cap B\bigl(0,K'\bigr) \bigr) \biggr| \biggr]\\
&&\qquad=0.
\end{eqnarray*}
Since, by choosing $K$ large enough, $P(\mathcal{V}^{[p]}\subset
B(0,p^{1/2}K))$ can be made arbitrarily close to $1$,
uniformly in $p$, we have proved that
%
%
%e32 #&#
\begin{eqnarray}
\label{BRW4} p^{-d/2} \#\mathcal{V}^{[p]} &\build{
\la}_{p\to\infty}^{\mathrm{(P)}}& \biggl(\frac{2}{\rho}
\biggr)^{d/2} \lambda_d \bigl(\bigl\{ M_\theta^{1/2}
\wh W_s \dvtx 0\leq s\leq\tau\bigr\} \bigr)
\nonumber
\\[-8pt]
\\[-8pt]
\nonumber
&= &\biggl(
\frac{2\sigma^2}{\rho} \biggr)^{d/2} \lambda_d \bigl( \{\wh
W_s \dvtx 0\leq s\leq\tau \} \bigr).
\end{eqnarray}
The relations between the Brownian snake and super-Brownian motion
\cite{Zurich}, Theorem IV.4, show that the quantity $\lambda_d(\{\wh W_s
\dvtx 0\leq s\leq\tau\})$
is the Lebesgue measure of the range of a super-Brownian motion (with
branching mechanism $2u^2$) started
from $(2/\rho)\delta_0$.
Finally, simple scaling arguments show that the limit can be expressed
in the form given in the theorem.
\end{pf}

%s5 #&#
\section{Open problems and questions}\label{sec5}

%s5.1 #&#
\subsection{The probability of visiting a distant point}\label{sec5.1}

Theorem~\ref{estim-visit} gives the asymptotic behavior of the
probability that a branching
random walk starting with a single particle at the origin visits a
distant point $a\in\Z^{d}$. It would
be of interest to have a similar result in dimension $d\geq4$,
assuming that $\theta$ is centered and has sufficiently high moments.
When $d\geq5$, a
simple calculation of the first and second moments of the number of
visits of $a$ (see, e.g., the remarks following Proposition~5 in \cite
{LGL}) gives the bounds
\[
C_1|a|^{2-d}\leq\Pi^*_{\mu,\theta}(a\in\bR)\leq
C_2 |a|^{2-d}
\]
with positive constants $C_1$ and $C_2$ depending on $d,\mu$ and
$\theta
$. When $d=4$,
one expects that
\[
\Pi^*_{\mu,\theta}(a\in\bR) \approx \frac{C}{|a|^2 \log|a|}.
\]
Calculations of moments give $\Pi^*_{\mu,\theta}(a\in\bR)\geq
c_1(|a|^2\log|a|)^{-1}$, but
proving the reverse bound $\Pi^*_{\mu,\theta}(a\in\bR)\leq
c_2(|a|^2\log|a|)^{-1}$
with some constant $c_2$ seems a nontrivial problem. This problem, in
the particular case of the
geometric offspring distribution, and some related questions are
discussed in Section~3.2 of \cite{BC}.

%s5.2 #&#
\subsection{The range in dimension four}\label{sec5.2}

With our previous notation $R_n$ for the range of a random walk indexed by
a random tree distributed according to $\Pi_\mu(\cdot \mid \#
\mathcal{T}
=n)$, Theorem~14 in \cite{LGL} states
that in dimension $d=4$,
\[
\frac{\log n}{n} R_n \build{\la}_{n\to\infty}^{L^2} 8
\pi^2\sigma^4,
\]
provided $\mu$ is the geometric distribution with parameter $1/2$, and
$\theta$ is symmetric and has exponential moments. It would be of
interest to extend this result
to more general offspring distributions. It seems difficult to adapt
the methods of
\cite{LGL} to a more general case, so new arguments would be needed.
In particular, finding the exact asymptotics of $\Pi^*_{\mu,\theta
}(a\in
\bR)$ (see
the previous subsection) in dimension $d=4$ would certainly be helpful.

%s5.3 #&#
\subsection{Branching random walk with a general initial
configuration}\label{sec5.3}
\label{bps}

One may ask whether a result such as Theorem~\ref{rangeBRW} remains
valid for
more general initial configurations of the branching particle system:
Compare with Propositions 20 and 21
in \cite{LGL}, which deal with the case $d\geq4$ and require no
assumption on the
initial configurations. In the present setting, Theorem~\ref{rangeBRW}
remains valid,
for instance, if we assume that the initial positions of the particles
stay within a bounded set
independently of $p$. On the other hand, one might consider the case
where we only assume that
the image of $p^{-1}X^{[p]}_0$ under the mapping $a\mapsto p^{-1/2}a$
converges weakly to a finite
measure $\xi$ on $\R^d$. This condition ensures the convergence of the
(rescaled) measure-valued processes
$X^{[p]}$ to a super-Brownian motion $Y$ with initial value $Y_0=\xi$,
and it is
natural to expect that we have, with a suitable constant~$C$,
%
%
%e33 #&#
\begin{equation}
\label{convbps} p^{-d/2} \#\mathcal{V}^{[p]} \build{
\la}_{p\to\infty}^{\mathrm{(d)}} C \lambda_d \biggl(\bigcup
_{t \geq0} \operatorname{ supp} Y_t \biggr).
\end{equation}

For trivial reasons, \eqref{convbps} will not hold in dimension $d=1$.
Indeed, for $\frac{1}{2}<\alpha<1$, we may let the initial
configuration consist of $p-\lfloor p^\alpha\rfloor$
particles uniformly spread over $\{1,2,\ldots,\sqrt{p}\}$ and
$\lfloor
p^\alpha\rfloor$ other particles located
at distinct points outside $\{1,2,\ldots,\sqrt{p}\}$. Then the
preceding assumptions hold ($\xi$ is the
Lebesgue measure on $[0,1]$), but \eqref{convbps} obviously fails since
$\#\mathcal{V}^{[p]}\geq\lfloor p^\alpha\rfloor$.
In dimension $2$, \eqref{convbps} fails again, for more subtle reasons:
One can construct examples where the descendants of certain initial
particles that
play no role in the convergence of the initial configurations
contribute to the asymptotics of $\#\mathcal{V}^{[p]}$
in a significant manner. Still, it seems likely that some version of
\eqref{convbps} holds under more stringent conditions on
the initial configurations [in dimension $3$ at least, the union in the
right-hand side of \eqref{convbps} should exclude $t=0$,
as can be seen from simple examples].

% imsref loaded by akundreckaite, 2014-07-15 12:15:41
%

%
%\begin{appendix}
%\section{}
%\end{appendix}

%\section*{Acknowledgments}

%\begin{supplement}[id=suppA]
%\slink[doi]{10.1214/00-AOPXXXXSUPP} %[doi,text={...}] - jei reikia
%suskaldyti doi
%\sfilename{aopXXXX\_supp.pdf}
%\end{supplement}

%\begin{thebibliography}{99}
%\bibitem[\protect\citeauthoryear{}{}]{r1}
%\bibitem{r1}
%\end{thebibliography}

\printaddresses
\end{document}